\documentclass[review,onefignum,onetabnum]{siamart190516}
\nolinenumbers

\usepackage{lipsum}
\usepackage{amsfonts}
\usepackage{graphicx}
\usepackage{epstopdf}
\usepackage{algorithmicx,dsfont} 
\usepackage{algpseudocode}
\usepackage{bm}
\ifpdf
  \DeclareGraphicsExtensions{.eps,.pdf,.png,.jpg}
\else
  \DeclareGraphicsExtensions{.eps}
\fi

\usepackage{cancel}
\usepackage{amsmath}
\usepackage{graphicx}
\usepackage[version=3]{mhchem}

\newsiamremark{remark}{Remark}
\newsiamremark{hypothesis}{Hypothesis}
\crefname{hypothesis}{Hypothesis}{Hypotheses}
\newsiamthm{claim}{Claim}

\headers{Subspace iteration with random sparsification}{S. M. Greene, R. J. Webber, T. C. Berkelbach, and J. Weare}

\title{Approximating matrix eigenvalues by subspace iteration with repeated random sparsification\thanks{ 
\funding{S.M.G. is supported by an investment fellowship from the Molecular Sciences Software Institute, which is funded by U.S. National Science Foundation grant OAC-1547580. R.J.W. is supported by New York University's Dean's Dissertation Fellowship and by the National Science Foundation through award DMS-1646339. J.W. acknowledges support from the Advanced Scientific Computing Research Program within the DOE Office of Science through award DE-SC0020427. The Flatiron Institute is a division of the Simons Foundation.}}}

\author{Samuel M. Greene\thanks{Department of Chemistry, Columbia University, New York, New York 10027, United States.}
\and Robert J. Webber\thanks{Courant Institute of Mathematical Sciences, New York University, New York, New York 10012, United States.}
  \and Timothy C. Berkelbach\footnotemark[2] \thanks{Center for Computational Quantum Physics, Flatiron Institute, New York, New York 10010, United States
  (\email{tim.berkelbach@gmail.com}).}
\and Jonathan Weare\footnotemark[3] \thanks{\email{weare@cims.nyu.edu}}}

\usepackage{amsopn}

\DeclareMathSymbol{\bot}{\mathord}{symbols}{"3F}

\ifpdf
\hypersetup{
  pdftitle={Subspace iteration with random sparsification},
  pdfauthor={S. M. Greene, R. J. Webber, T. C. Berkelbach, and J. Weare}
}
\fi

\begin{document}

\maketitle

\begin{abstract}
  Traditional numerical methods for calculating matrix eigenvalues are prohibitively expensive for high-dimensional problems.
Iterative random sparsification methods allow for the estimation of a single dominant eigenvalue at reduced cost by leveraging repeated random sampling and averaging.
We present a general approach to extending such methods for the estimation of multiple eigenvalues and demonstrate its performance for several benchmark problems in quantum chemistry.
\end{abstract}

\begin{keywords}
  eigenvalues, subspace iteration, randomized algorithms, Monte Carlo
\end{keywords}

\begin{AMS}
  65F15, 68W20, 65C05
\end{AMS}

\section{Introduction}
A wide range of applications, including principal component analysis \cite{Wold1987, Herve2010principal}, spectral analysis of dynamical systems \cite{Noe2013, williams2015data}, and electronic structure calculations \cite{Knowles1984, Bickelhaupt2000}, 
require matrix eigenvectors and eigenvalues. Methods for calculating them based on dense, in-place factorizations are intractably expensive for large matrices \cite{Stewart1998vol1, Stewart2001vol2}.
Iterative methods involving repeated matrix--vector multiplications \cite{Lanczos1950, Davidson1975, Sleijpen1996} offer reduced computational and memory costs, particularly for sparse matrices.
However, even these methods are too expensive for the extremely large matrices increasingly encountered in modern applications. 

We consider a class of randomized iterative techniques that enable significant further reductions in memory and computational costs.
These techniques build on classical iterative methods for solving eigenvalue problems 
by randomly perturbing vectors or matrices at each iteration to increase sparsity.
The imposed sparsity facilitates the use of sparse linear algebra frameworks for performing matrix--vector multiplication efficiently. Because they do not involve storing or manipulating dense vectors, iterative random sparsification methods
are particularly suited to problems so large that storing even a single dense vector is unmanageably expensive \cite{Lim2017}.
Iterative random sparsification methods can estimate the dominant eigenvalues and eigenspaces of extremely large matrices at reduced cost \cite{Mascagni2013rnla,Lim2017,Lu2017,Andoni2018solving,ozdaglar2019asynchronous} compared to randomized or deterministic methods that repeatedly apply matrix--vector multiplications with dense vectors \cite{rokhlin2010randomized,halko2011algorithm,halko2011finding,gu2015subspace,musco2015randomized}.

When used to estimate the ground-state energy, or smallest eigenvalue, of the quantum mechanical Hamiltonian operator, iterative random sparsification methods are termed ``projector quantum Monte Carlo'' methods~\cite{Nightingale1996, Cleland2010, Lu2017, Motta2018, Greene2019, Greene2020}.
These methods have become a standard tool for calculating ground-state electronic energies of molecules.
Applying such methods to calculate \emph{multiple} eigenvalues poses additional challenges related to the need to maintain orthogonality among eigenvectors as the iteration proceeds \cite{Ceperley1988, Ohtsuka2010, Blunt2015, Feldt2020}. 
Yet successful extension to the multiple eigenvalue problem would have implications not only for quantum applications, but in other areas of physics, engineering, and data science.

This article presents a subspace iteration approach with repeated random sparsification for the estimation of multiple dominant eigenvalues. 
The approach is based on a version of subspace iteration with several non-standard design choices to increase the method's stability under random perturbations. 
Random perturbations are introduced to the vectors at each iteration to promote sparsity. These perturbations have the smallest possible mean square magnitude,
while preserving the original vector in expectation.
The resulting randomized subspace iteration builds on previous methods for calculating the single dominant eigenvalue within the fast randomized iteration framework~\cite{Lim2017, Lu2017, Greene2019, Greene2020}, but the extension from estimating one dominant eigenvalue to multiple dominant eigenvalues is new, as is the random sparsification technique.

The randomized subspace iteration introduced here has traits in common with Markov chain Monte Carlo \cite{Levin2017mixing}.
First, as discussed  in Section \ref{sec:theory}, the matrix should be \emph{irreducible} to ensure the method's stability.
Second, the random iterates generated by the scheme converge to a stationary distribution, rather than the desired (deterministic) eigenspace itself. 
Statistical error in the eigenvalue estimates can be reduced by averaging over many iterations $i$, resulting in a convergence rate of $1 \slash \sqrt{i}$.
However, the estimates are not asymptotically consistent as $i \rightarrow \infty$.
Rather, eigenvalue estimates include an asymptotic bias that can be systematically reduced by relaxing the sparsity constraint.

We test our method on the full configuration interaction eigenproblem from many-electron quantum mechanics; in this context, it can be understood as a generalization of projector quantum Monte Carlo methods to excited states.
Our tests indicate that the new algorithm reduces per-iteration computational costs compared to deterministic subspace iteration
and yields eigenvalue estimates of sufficient accuracy for the quantum mechanical problems considered here ($<1$ m$E_\text{h}$ errors).
This high accuracy is especially remarkable because the number of nonzero entries retained at each iteration is less than $1\%$ of the dimension of the matrix, even for these small examples.
In a recent companion work \cite{greene2022full}, we have extended the method to obtain reliable eigenvalue estimates for dramatically larger matrices, with up to $10^{25} \times 10^{25}$ entries.

The remainder of this paper begins with a description of the deterministic subspace iteration method on which our randomized algorithm is based, with a particular focus on its non-standard features that make it robust to random perturbations (Section \ref{sec:detSubsp}).
Section \ref{sec:comp} describes our approach to stochastically imposing sparsity in vectors and matrices.
Section \ref{sec:randSubsp} describes our randomized subspace iteration scheme.
Section \ref{sec:theory} presents a theoretical error analysis.
Section \ref{sec:results} presents applications of this algorithm to quantum mechanical problems.
In Section \ref{sec:conclusions}, we summarize our key findings and discuss possible methodological improvements.

Throughout this work, we use the following notation:
\begin{itemize}
\item Matrices $\textbf{X}$ are written in bold capital letters, vectors $\textbf{x}$ are written in bold lower case letters,
and scalars $x$ are written in italic lower case letters.
\item We use $\textbf{x}_i$ to indicate the $i^\text{th}$ entry of a vector,
$\textbf{X}_{ij}$ to indicate the $\left(i,j\right)\text{th}$ entry of a matrix, and $\textbf{X}_{:i}$ to indicate the $i$th column of a matrix.
\item For any vector $\mathbf{x} \in \mathbb{R}^n$, $\left\lVert \mathbf{x} \right\rVert_2 = (\sum_{i=1}^n |\mathbf{x}_i|^2)^{1 \slash 2}$ denotes the Euclidean norm,
$\left\lVert \mathbf{x} \right\rVert_1 = \sum_{i=1}^n |\mathbf{x}_i|$
denotes the sum of the absolute values of the entries,
and $\left\lVert \mathbf{x} \right\rVert_0 = \sum_{i=1}^n \mathds{1}\left\{|\mathbf{x}_i \neq 0\right\}$ denotes the number of nonzero entries.
\item For any matrix $\mathbf{X} \in \mathbb{R}^{m \times n}$, $\lVert \mathbf{X} \rVert_{\textup{F}} = ( \sum_{i=1}^m \sum_{j=1}^n \mathbf{X}_{ij}^2 )^{1 \slash 2}$ denotes the Frobenius norm.
\item $\mathbf{X}|_S$ signifies the restriction of a matrix $\mathbf{X}$ to a linear subspace $S$, while $\mathbf{P}_S$ and $\mathbf{P}_{\mathbf{X}}$ signify orthogonal projections onto a linear subspace $S$ or the range of a matrix $\mathbf{X}$.
\item $\textup{E}$, $\textup{P}$, $\textup{Var}$, and $\textup{Cov}$ indicate expectations, probabilities, variances, and covariances with respect to a probability space that is sufficiently rich to support all the random variables identified in the analysis.
\end{itemize}

\section{A non-standard deterministic subspace iteration}
\label{sec:detSubsp}
Our goal is to find the $k$ dominant eigenvalues (counting multiplicity) of a matrix $\mathbf{A} \in \mathbb{R}^{n \times n}$. Starting from an initial matrix $\mathbf{X}^{(0)} \in \mathbb{R}^{n\times k}$, standard subspace iteration constructs a sequence of matrix iterates according to the recursion 
\begin{equation}
\label{eq:detIteration}
\mathbf{X}^{(i + 1)} = \mathbf{A} \mathbf{X}^{(i)} [\mathbf{G}^{(i)}]^{-1},
\end{equation}
where multiplication by $[\mathbf{G}^{(i)}]^{-1}$ enforces orthonormality among the columns of $\mathbf{X}^{(i + 1)}$ \cite{Stewart1969, Stewart1975, Saad2011}.
For $k=1$, subspace iteration reduces to power iteration, on which many single-eigenvalue randomized iterative methods are based.
For $k > 1$ (the case considered here), subspace iteration provides a solution to the multiple dominant eigenvalue problem.
Eigenvalues of $\mathbf{A}$ can be estimated after each iteration
by solving the eigenvalue problem
\begin{equation}
\label{eq:symmRayRitz}
\mathbf{X}^{(i)*} \mathbf{A} \mathbf{X}^{(i)} \mathbf{W}^{(i)} = \mathbf{W}^{(i)} \mathbf{\Lambda}^{(i)}
\end{equation}
for the diagonal matrix $\mathbf{\Lambda}^{(i)}$ of eigenvalue estimates~\cite{Stewart1969}.

Standard subspace iteration involves nonlinear operations on the iterates, both for enforcing the orthonormality of columns 
and for estimating eigenvalues by \eqref{eq:symmRayRitz}.
However, these operations lead to statistical biases once randomness is introduced into the iterates by stochastic sparsification.
In order to reduce these errors in our randomized algorithm, we make two \emph{non-standard} choices.

As a first non-standard choice, we estimate eigenvalues by solving the eigenvalue problem
\begin{equation}
\label{eq:genEval}
\mathbf{U}^{*} \mathbf{A} \mathbf{X}^{(i)} \mathbf{W}^{(i)} = \mathbf{U}^{*} \mathbf{X}^{(i)} \mathbf{W}^{(i)} \mathbf{\Lambda}^{(i)},
\end{equation}
for the diagonal matrix $\mathbf{\Lambda}^{(i)}$ of eigenvalue estimates, 
where $\mathbf{U}$ is a constant deterministic matrix with columns chosen to approximate the dominant eigenvectors of $\mathbf{A}$.
Equation (\ref{eq:genEval}) generalizes the ``projected estimator'' commonly used in single-eigenvalue randomized methods~\cite{Booth2014}.
Here, we also use  $\mathbf{U}$ as the first iterate, i.e.~$\mathbf{X}^{(0)} = \mathbf{U}$.

The matrix $\mathbf{U}$ must be obtained by other means before application of subspace iteration.
On the one hand,  $\mathbf{U}$ should be chosen to approximate the subspace of the dominant $k$ eigenvectors as nearly as possible,
both to ensure a reasonable starting iterate and to optimize eigenvalue accuracy at every estimation step.
Indeed, eigenvalue estimates are exact for any eigenvector contained within the column span of $\mathbf{U}$.
On the other hand, when $n$ is extremely large, assembling the matrices in \eqref{eq:genEval} may require imposing sparsity or some other structural restriction on $\mathbf{U}$.
Explicit error bounds that reveal the importance of choosing a good starting matrix $\mathbf{U}$ are stated and proved in Section \ref{sec:deterministicAnalysis}.


As a second non-standard choice,
at intervals of $\Delta$ iterations,
we construct the matrices $\mathbf{G}^{(i)}$
so that multiplication by $[\mathbf{G}^{(i)}]^{-1}$ makes the columns of the projected matrix $\mathbf{P}_{\mathbf{U}} \mathbf{X}^{(i+1)}$ orthogonal, rather than enforcing the orthogonality of $\mathbf{X}^{(i+1)}$.
Specifically, we set 
\begin{equation}
\mathbf{G}^{(i)} = \mathbf{N}^{(i)} \mathbf{D}^{(i)} \mathbf{R}^{(i)}
\end{equation}
where $\mathbf{R}^{(i)}$ is the upper triangular factor from a QR factorization of $\mathbf{U}^* \mathbf{A} \mathbf{X}^{(i)}$,
and $\textbf{D}^{(i)}$ is the diagonal matrix with entries
$\textbf{D}^{(i)}_{jj} = 
\bigl\lVert (\mathbf{X}^{(i)} [\mathbf{R}^{(i)}]^{-1})_{:j} \bigr\rVert_1 \slash \bigl\lVert \mathbf{X}^{(i)}_{:j} \bigr\rVert_1$.
The diagonal matrix  $\mathbf{N}^{(i)}$, defined recursively by
\begin{equation}
\textbf{N}^{(i)}_{jj} = \Biggl( \frac{\bigl\lVert \mathbf{X}^{(i)}_{:j} \bigr\rVert_1}
{\bigl\lVert \mathbf{X}^{(i-1)}_{:j} \bigr\rVert_1} \Biggr)^\alpha \bigl( \textbf{N}^{(i-1)}_{jj} \bigr)^{(1 - \alpha)}, \quad \textbf{N}^{(0)}_{jj} = 1
\end{equation}
controls the scaling of the columns of $\mathbf{X}^{(i+1)}$. 
Choosing the user-defined parameter $\alpha$ equal to 1 keeps the $\ell_1$-norm of each iterate column constant (i.e.~$\bigl\lVert \mathbf{X}^{(i)}_{:j} \bigr\rVert_1 = \bigl\lVert \mathbf{X}^{(i - 1)}_{:j} \bigr\rVert_1$) but introduces a bias since $\mathbf{N}^{(i)}$ depends non-linearly on iterates.
Choosing $\alpha < 1$ reduces bias by damping this non-linear dependence.
We found that $\alpha = 0.5$ was a suitable choice for our numerical experiments.
For iterations $i$ at which we do not orthogonalize, we still control the scaling of the columns of $\mathbf{X}^{(i+1)}$ by setting
$\mathbf{G}^{(i)} = \mathbf{N}^{(i)}$.


The interval $\Delta$ at which orthogonalization should be performed can be determined by monitoring the condition number of the matrix $\mathbf{U}^* \mathbf{X}^{(i)}$.
If orthogonalization is performed too infrequently, the condition number will increase as the iteration proceeds, giving rise to instabilities in the  algorithm.
For the specific numerical experiments in this article, we monitored this condition number (Appendix \ref{sec:stability}) and found $\Delta = 1000$ iterations to be a suitable orthogonalization interval. Our results did not change significantly when we orthogonalized more frequently.
However, orthogonalization will need to be performed more frequently --- perhaps as often as every iteration ($\Delta = 1$) --- for matrices with faster-decaying spectra than those in our experiments.

The non-standard eigenvalue estimator and orthogonalization step would lead to a suboptimal deterministic algorithm.
Indeed, if $\mathbf{A}$ is symmetric and positive definite, eigenvalues estimated using our non-standard approach have errors that decays at a rate $( \lambda_j / \lambda_{k+1})^{i}$ as $i \rightarrow \infty$, whereas the standard eigenvalue estimator \eqref{eq:symmRayRitz} leads to errors that decay faster, at a rate $( \lambda_j / \lambda_{k+1})^{2i}$ (Section \ref{sec:deterministicAnalysis}).
However, when random perturbations are applied to the matrix iterates, our choices lead to a more stable algorithm.
The standard subspace iteration is based on quadratic forms
$\mathbf{X}^{(i)*} \mathbf{A X}^{(i)}$ and $\mathbf{X}^{(i)*} \mathbf{X}^{(i)}$ that
have a high bias and variance in the presence of stochastic perturbations \cite{Lim2017}.
In contrast, our non-standard subspace iteration is based on quantities that are better-behaved, including the linear products of random matrices and constant matrices  $\mathbf{U}^* \mathbf{X}^{(i)}$
and $\mathbf{U}^* \mathbf{A} \mathbf{X}^{(i)}$, as well as the random column norms $||\mathbf{X}_{:j}^{(i)}||$.
As a result, our randomized subspace iteration can lead to highly accurate eigenvalue estimates (Section \ref{sec:results}), whereas iterative random sparsification applied to the standard subspace iteration fails dramatically (Appendix \ref{sec:symmSubsp}).


\section{Stochastic compression}
\label{sec:comp}
If sparsity is leveraged, the cost of forming the matrix products $\mathbf{A X}^{(i)}$ in the above algorithm scales as $\mathcal{O}(m_a m_x k)$, where $m_a$ and $m_x$ are the maximum number of nonzero elements in each column of $\mathbf{A}$ and $\mathbf{X}^{(i)}$, respectively. 
Stochastic compression allows one to control this cost by zeroing nonzero elements at randomly selected positions.
We define a stochastic compression operator $\Phi$ which, when applied to a vector $\mathbf{x}$, returns a random compressed vector $\Phi(\mathbf{x})$ with (1) at most a user-specified number $m$ of nonzero elements and (2) all elements equal to those of the input vector $\mathbf{x}$ in expectation, i.e.,  E$[\Phi(\mathbf{x})_i] = \mathbf{x}_i$.
Applying $\Phi$ to a matrix $\mathbf{X} = [\mathbf{X}_{:1} \quad \mathbf{X}_{:2} \quad ...]$ involves compressing each of its columns independently to $m$ nonzero elements, i.e.~$\Phi(\mathbf{X}) = [ \Phi(\mathbf{X}_{:1}) \quad \Phi(\mathbf{X}_{:2}) \quad ...]$.

In the remainder of this section, we describe a compression scheme called ``pivotal compression''.
We find that pivotal compression achieves less variance in applications than previous compression schemes~\cite{Lim2017, Greene2019, Greene2020}, and we prove in Section \ref{sec:compressionAnalysis} that
this scheme minimizes the mean square perturbation magnitude $\textup{E}\left\lVert \Phi(\mathbf{x}) - \mathbf{x} \right\rVert_2^2$
over all possible compression schemes.
Additionally, when $m$ equals or exceeds the number of nonzero elements in $\mathbf{x}$, pivotal compression yields the exact input vector $\mathbf{x}$, in which case the statistical error is zero.

The first step in pivotal compression involves identifying a set $\mathcal{D}$ of indices corresponding to the $d$ largest-magnitude elements in $\mathbf{x}$.
These elements are preserved exactly during the compression, i.e.~$\Phi(\mathbf{x})_i = \mathbf{x}_i$ for all $i \in \mathcal{D}$. 
The number $d$ is chosen such that
\begin{align}
\label{eq:detCriterion}
& |\mathbf{x}_i| \geq \frac{1}{m-d} \sum_{j \notin \mathcal{D}} |\mathbf{x}_j|, \qquad \forall i \in \mathcal{D}, \\
\label{eq:detCriterion2}
& |\mathbf{x}_i| \leq \frac{1}{m-d} \sum_{j \notin \mathcal{D}} |\mathbf{x}_j|, \qquad \forall i \notin \mathcal{D}.
\end{align}
The next step is to 
apply pivotal sampling~\cite{deville1998unequal, chauvet2012characterization, Chauvet2017} to
randomly select a set $\mathcal{S}$ consisting of $m - d$ additional indices.
The probability that each index $i$ is included in $\mathcal{S}$ is 
\begin{equation}
\label{eq:probabilities}
\mathbf{p}_i = 
\frac{\left(m - d\right) |\mathbf{x}_i|}{\sum_{j \notin \mathcal{D}} |\mathbf{x}_j|}, \qquad \forall i \notin \mathcal{D}.
\end{equation}
Lastly, the sets $\mathcal{D}$ and $\mathcal{S}$ are used to construct the compressed vector $\Phi\left(\mathbf{x}\right)$ as
\begin{equation}
\Phi(\mathbf{x})_i = \begin{cases}
\mathbf{x}_i, & i \in \mathcal{D}, \\
\mathbf{x}_i \slash \mathbf{p}_i, & i \in \mathcal{S}, \\
0, & i \notin \mathcal{D}, i \notin \mathcal{S},
\end{cases} 
\end{equation}
thus ensuring that this compression scheme is unbiased (i.e.~that $\textup{E} \Phi(\mathbf{x}) = \mathbf{x}$).
We provide complete pseudocode for performing pivotal compression in Appendix \ref{sec:pivComp} and discuss a possible approach to parallelizing pivotal compression in Appendix \ref{sec:parComp}.

\section{Randomized subspace iteration}
\label{sec:randSubsp}
Combining stochastic compression with our subspace iteration yields the randomized subspace iteration
\begin{equation}
\label{eq:iteration}
    \mathbf{X}^{(i + 1)} = \mathbf{A} \Phi(\mathbf{X}^{(i)}) [\mathbf{G}^{(i)}]^{-1},
\end{equation} 
where the compression operation $\Phi$ is performed independently at each iteration. 
The matrix products $\mathbf{U}^* \mathbf{A} \Phi(\mathbf{X}^{(i)})$ and $\mathbf{U}^* \mathbf{X}^{(i)}$ are evaluated and stored for the purpose of estimating eigenvalues.
We provide a complete pseudocode in Algorithm \ref{alg:randSubsp} below.

\begin{algorithm}
\caption{Randomized subspace iteration}
\label{alg:randSubsp}
\begin{algorithmic}[1]
\State{\textbf{Input:} Matrix $\mathbf{U}$ whose columns approximate the $k$ dominant eigenvectors of $\mathbf{A}$}
\State{\textbf{Initialization:} $\mathbf{X}^{(0)} = \mathbf{U}$ and $\mathbf{N}^{(0)} = \mathbf{I}$}
\For{$i=0, 1, 2, ..., i_\text{max} - 1$} 
\State{$\mathbf{J}^{(i)} = \mathbf{U}^* \mathbf{X}^{(i)}$}
\State{Construct $\mathbf{X}^{(i)\prime} = \Phi(\mathbf{X}^{(i)})$} \Comment{Pivotal compression}
\State{$\mathbf{Y}^{(i)} = \mathbf{A} \mathbf{X}^{(i)\prime}$}
\State{$\mathbf{K}^{(i)} = \mathbf{U}^* \mathbf{Y}^{(i)}$}
\If{$i \equiv (\Delta - 1)\, (\text{mod } \Delta)$}
\State{Construct $\mathbf{Q}^{(i)} \mathbf{R}^{(i)} = \mathbf{K}^{(i)}$} \Comment{QR factorization}
\State{$\mathbf{Z}^{(i)} = \mathbf{Y}^{(i)} [\mathbf{R}^{(i)}]^{-1}$}
\For{$j=1, 2, 3, ..., k$}
\State{$\mathbf{D}^{(i)}_{jj} = \bigl\lVert \mathbf{Z}_{:j}^{(i)} \bigr\rVert_1 \bigl\lVert \mathbf{Y}_{:j}^{(i)}\bigr\rVert_1^{-1}$} 
\EndFor
\State{$\mathbf{G}^{(i)} = \mathbf{N}^{(i)} \mathbf{D}^{(i)} \mathbf{R}^{(i)}$}
\Else
\State{$\mathbf{G}^{(i)} = \mathbf{N}^{(i)}$}
\EndIf
\State{$\mathbf{X}^{(i + 1)} = \mathbf{Y}^{(i)} [\mathbf{G}^{(i)}]^{-1}$}
\For{$j = 1, 2, ..., k$}
\State{$\mathbf{N}^{(i+1)}_{jj} = 
\bigl\lVert \mathbf{X}^{(i+1)}_{:j} \bigr\rVert_1^{\alpha} 
\bigl\lVert \mathbf{X}^{(i)}_{:j} \bigr\rVert_1^{-\alpha} 
\bigl[\textbf{N}^{(i-1)}_{jj}\bigr]^{1 - \alpha}$.}
\EndFor
\EndFor
\State{\textbf{Return:} Sequence of matrices $\mathbf{J}^{(i)} = \mathbf{U}^* \mathbf{X}^{(i)}$ and $\mathbf{K}^{(i)} = \mathbf{U}^* \mathbf{A} \Phi(\mathbf{X}^{(i)})$.}
\end{algorithmic}
\end{algorithm}

In order to produce eigenvalue estimates
after running Algorithm \ref{alg:randSubsp}, we first form the matrices
\begin{align}
\label{eq:ave1}
    &\langle \mathbf{U}^* \mathbf{A} \Phi(\mathbf{X}) \rangle_{i_{\max}} = \frac{1}{i_\text{max} - i_\text{min}} \sum_{i = i_\text{min}}^{i_\text{max} - 1} \mathbf{U}^* \mathbf{A} \Phi(\mathbf{X}^{(i)}), \\
\label{eq:ave2}
    & \langle \mathbf{U}^* \mathbf{X} \rangle_{i_{\max}}
    = \frac{1}{i_\text{max} - i_\text{min}} \sum_{i = i_\text{min}}^{i_\text{max} - 1} \mathbf{U}^* \mathbf{X}^{(i)},
\end{align}
which are time-averages starting from an initial burn-in time $i_\text{min}$ and going until a final time $i_\text{max}$.
Then, we solve the eigenvalue equation
\begin{equation}
\label{eq:linEst}
\langle \mathbf{U}^* \mathbf{A} \Phi(\mathbf{X}) \rangle_{i_{\max}} \mathbf{W} = \langle \mathbf{U}^* \mathbf{X} \rangle_{i_{\max}} \mathbf{W} \mathbf{\Lambda}^{(i_\text{max})}.
\end{equation}


From our theoretical analysis in Section \ref{sec:randSubspAnalysis}, we anticipate that the time averages \eqref{eq:ave1} and \eqref{eq:ave2} will converge as $i_{\max} \rightarrow \infty$, for any fixed value of the burn-in time $i_{\min} \geq 0$.
However, 
it is beneficial
to choose a positive burn-in time $i_{\min} > 0$, which reduces the period of \emph{initialization bias} in which \eqref{eq:ave1} and \eqref{eq:ave2} systematically deviate from their limiting values \cite{Sokal1997}.
To choose $i_{\min}$ in our numerical experiments,
we use the approach of \cite{Chodera2016}, whereby we increase $i_{\min}$ until the autocorrelation times associated with our eigenvalue estimates begin to stabilize
(signalling that the initialization period is sufficiently long).

Lastly, the analysis in Section \ref{sec:randSubspAnalysis} indicates that the standard error of each eigenvalue estimate $\mathbf{\Lambda}_{jj}^{(i_\text{max})}$ can be estimated through the formula
\begin{equation}
\label{eq:var_calc}
\textup{Var}[\mathbf{\Lambda}_{jj}^{(i_\text{max})}]
\approx \frac{1}{(i_\text{max} - i_\text{min})^2}
\sum_{\substack{i_{\min} \leq i,k \leq i_{\max} - 1 \\
\left|i - k\right| \leq \tau
}}
f_j^{(i)} f_j^{(k)}.
\end{equation}
Here, $f^{(i_{\min})}_j, f^{(i_{\min + 1})}_j, \ldots$ is a scalar-valued time series, defined by
\begin{equation}
f^{(i)}_j = {\mathbf{z}}_j^* \left(\mathbf{U}^* \mathbf{Y}^{(i)} - \mathbf{\Lambda}_{jj}^{(i_\text{max})} \mathbf{U}^* \mathbf{X}^{(i)} \right) {\mathbf{w}}_j,
\end{equation}
and $\mathbf{z}_k$ and $\mathbf{w}_k$ represent left and right generalized eigenvectors corresponding to $\mathbf{\Lambda}_{jj}^{(i_\text{max})}$.
The truncation threshold $\tau > 0$ is chosen so that correlations involving $f_j^{(i)}$ and $f_j^{(k)}$ are negligibly small for any $|i - k| > \tau$ \cite{Sokal1997}.
Determining the truncation threshold $\tau$ and computing the variance \eqref{eq:var_calc} are common procedures in Markov chain Monte Carlo and are conveniently implemented in the emcee package for python \cite{Foreman2013},
which we use throughout our numerical experiments.


\section{Theoretical analysis}\label{sec:theory}

In this section, we mathematically justify several of the design components of randomized subspace iteration (Algorithm \ref{alg:randSubsp}).
First, we prove that Algorithm \ref{alg:randSubsp} is based on a deterministic subspace iteration that converges exponentially fast (Section \ref{sec:deterministicAnalysis}).
Then, we prove that the random perturbations are as small as possible, in the sense of minimizing $\textup{E} | \mathbf{\Phi}(\mathbf{x}) - \mathbf{x} |^2$ over all possible compression schemes (Section \ref{sec:compressionAnalysis}).
Lastly, we derive an \emph{a posteriori} variance estimator for the estimated eigenvalues, assuming geometric ergodicity (Section \ref{sec:randSubspAnalysis}).

These results confirm that Algorithm \ref{alg:randSubsp} is built on rigorous principles, but they do not provide a complete convergence theory for the algorithm.
In the future, more analysis is needed to precisely identify the set of problems for which randomized subspace iteration performs well,
as well as identifying  failure modes for the new method.

One failure mode for randomized subspace iteration is already known.
When the approach is applied to a block diagonal matrix, the elements of $\mathbf{X}^{(i)}$ in a single block can become zero,
in which case they remain zero for all time.
This makes Algorithm \ref{alg:randSubsp} less stable than the corresponding deterministic scheme.
In a previous analysis, Lim \& Weare \cite{Lim2017} avoid the degeneracies associated with block diagonal matrices by focusing
on \emph{irreducible} matrices $\mathbf{A}$, i.e., matrices for which every basis element is connected to every other basis element by a chain of nonzero matrix entries.
Assuming irreducibility (and aperiodicity), \cite{Lim2017} establishes an error bound for a simplified Algorithm \ref{alg:randSubsp} that estimates the single dominant eigenvalue of a matrix with nonnegative entries.
However, it remains uncertain whether irreducibility is enough to guarantee stability in the general case, which may involve multiple eigenvalues and matrices $\mathbf{A}$ with both positive and negative entries.

The remainder of the section presents mathematical arguments supporting the design components of Algorithm \ref{alg:randSubsp}.
These arguments are direct but slightly lengthy, so the reader may want to skip over the proofs on a first reading.
We defer to Section \ref{sec:results} the numerical experiments that show our approach is effective on realistic problems.



\subsection{Analysis of the deterministic subspace iteration}
\label{sec:deterministicAnalysis}

Here we derive error bounds for the non-standard deterministic subspace iteration presented in Section \ref{sec:detSubsp}, assuming the matrix $\mathbf{A} \in \mathbb{R}^{n \times n}$ is symmetric.
These error bounds are based on two observations.

The first observation is that the iteration 
\eqref{eq:detIteration} leads to an explicit representation for the matrix $\mathbf{X}^{(i)}$ as
\begin{equation}
\label{eq:recursive}
\mathbf{X}^{(i)} = \mathbf{A}^i \mathbf{U} \left[\mathbf{G}^{(i-1)} \mathbf{G}^{(i-2)} \cdots \mathbf{G}^{(0)}\right]^{-1}.
\end{equation}
Therefore, by substituting $\mathbf{V}^{(i)} = \left[\mathbf{G}^{(i-1)} \mathbf{G}^{(i-2)} \cdots \mathbf{G}^{(0)}\right]^{-1} \mathbf{W}^{(i)}$
into \eqref{eq:genEval}, the eigenvalue problem simplifies to become
\begin{equation}
\label{eq:symmetric}
    \underbrace{\mathbf{U}^{\ast} \mathbf{A}^{i+1} \mathbf{U}}_{\text{symmetric matrix}} \mathbf{V}^{(i)} = \underbrace{\mathbf{U}^{\ast} \mathbf{A}^i \mathbf{U}}_{\text{symmetric matrix}} \mathbf{V}^{(i)} \mathbf{\Lambda}^{(i)},
\end{equation}
which is a symmetric generalized eigenvalue problem involving matrices $\mathbf{U}^{\ast} \mathbf{A}^i \mathbf{U}$ and $\mathbf{U}^{\ast} \mathbf{A}^{i+1} \mathbf{U}$.

The second observation is that the eigenvalues estimates in \eqref{eq:symmetric} satisfy the Courant--Fischer min--max principle \cite{reed1979methods,parlett1998symmetric},
provided that $\mathbf{A}^i$ is positive definite and $\mathbf{U}$ is nonsingular.
Indeed, if eigenvalue estimates are ordered 
from largest to smallest,
\begin{equation}
    \mathbf{\Lambda}_{11}^{(i)} \geq \mathbf{\Lambda}_{22}^{(i)} \geq \cdots \geq \mathbf{\Lambda}_{kk}^{(i)},
\end{equation}
the min--max principle gives the representation
\begin{equation}
\label{eq:variational}
    \mathbf{\Lambda}_{jj}^{(i)} = \max_{\substack{S \subseteq \textup{range}[\mathbf{U}] \\ \textup{dim}[S] = j}}
    \min_{\mathbf{x} \in S} \frac{\mathbf{x}^{\ast} \mathbf{A}^{i+1} \mathbf{x}}{\mathbf{x}^{\ast} \mathbf{A}^{i} \mathbf{x}}, \qquad 1 \leq j \leq k.
\end{equation}

The min--max representation \eqref{eq:variational}
is similar to a variational representation for the true eigenvalues of $\mathbf{A}$, but the only difference is the dependence on $\textup{range}[\mathbf{U}]$.
Indeed, if we order the eigenvalues of $\mathbf{A}$ from largest to smallest,
\begin{equation}
    \lambda_1 \geq \lambda_2 \geq \cdots \geq \lambda_n,
\end{equation}
the eigenvalues of $\mathbf{A}$ satisfy
\begin{equation}
\label{eq:variational2}
    \lambda_j = \max_{\substack{S \subseteq \mathbb{R}^n \\ \textup{dim}[S] = j}}
    \min_{\mathbf{x} \in S} \frac{\mathbf{x}^{\ast} \mathbf{A}^{i+1} \mathbf{x}}{\mathbf{x}^{\ast} \mathbf{A}^{i} \mathbf{x}}, \qquad 1 \leq j \leq n,
\end{equation}
which is an alteration of \eqref{eq:variational} with $\textup{range}[\mathbf{U}]$ replaced by the larger space $\mathbb{R}^n$. 
Thus, by comparing \eqref{eq:variational} and \eqref{eq:variational2} we obtain the inequality
\begin{equation}
    \mathbf{\Lambda}^{(i)}_{jj} \leq \lambda_j,
\end{equation}
valid for all $1 \leq j \leq k$.

So far, we have bounded the eigenvalue estimates $\mathbf{\Lambda}^{(i)}_{jj}$ from above using the true eigenvalues $\lambda_j$.
However, the task remains to bound the eigenvalues $\mathbf{\Lambda}^{(i)}_{jj}$ from below, which is exactly what we accomplish in the next proposition.

\begin{proposition}
\label{prop:improved}
Let $\beta$ be a permutation that reorders the eigenvalues of $\mathbf{A}$ by magnitude, i.e.,
\begin{equation}
    \left|\lambda_{\beta(1)}\right| \geq \left|\lambda_{\beta(2)}\right| \geq \cdots \geq \left|\lambda_{\beta(n)}\right|,
\end{equation}
and let $H_k$ be a subspace of eigenvectors of $\mathbf{A}$ with eigenvalues $\lambda_{\beta(1)}, \ldots, \lambda_{\beta(k)}$.
Let
\begin{equation}
    \label{eq:theta_def}
    \theta = \arccos\left( \min_{\mathbf{x} \in H_k} 
    \frac{\left\lVert \mathbf{P}_{\mathbf{U}} \mathbf{x}\right\rVert}{\left\lVert \mathbf{x} \right\rVert}\right)
\end{equation}
be the angular distance between $H_k$ and $\textup{range}[\mathbf{U}]$,
and assume $\theta < \frac{\pi}{2}$.
Lastly, assume that $\mathbf{A}^i$ is positive definite and $\lambda_j \geq 0$.
Then, the eigenvalue estimate $\mathbf{\Lambda}^{(i)}_{jj}$ resulting from \eqref{eq:detIteration} and \eqref{eq:genEval}
satisfies
\begin{equation}
\label{eq:explicit}
    \frac{1 - R^{i+1} \tan^2 \theta}
    {1 + R^i \tan^2 \theta} \lambda_j \leq \mathbf{\Lambda}^{(i)}_{jj} \leq \lambda_j \quad \text{where} \quad R = \left|\frac{\lambda_{\beta(k+1)}}{\lambda_j}\right|.
\end{equation}
\end{proposition}
\begin{proof}
The idea of the proof is to construct a particular subset $S \subseteq \textup{range}[\mathbf{U}]$ such that $\textup{dim}[S] \geq j$ and
\begin{equation}
\label{eq:verify}
    \frac{\mathbf{x}^{\ast} \mathbf{A}^{i+1} \mathbf{x}}{\mathbf{x}^{\ast} \mathbf{A}^{i} \mathbf{x}} \geq \frac{1 - \tan^2 \theta R^{i+1}}
    {1 + \tan^2 \theta R^i} \lambda_j, \qquad \forall x \in S.
\end{equation}
Then the Courant--Fischer min--max principle \eqref{eq:variational} guarantees \eqref{eq:explicit}.

To build the set $S$, let
$J \subseteq H_k$ be the subspace spanned by eigenvectors of $\mathbf{A}$ with eigenvalues $\lambda_1, \ldots, \lambda_j$,
and define
\begin{equation}
    J^{\prime} = H_k \cap J^{\bot},
    \qquad Q_k = \mathbf{P}_{\mathbf{U}} H_k,
    \qquad S = Q_k \cap (\mathbf{P}_{\mathbf{U}} J^{\prime})^{\bot}.
\end{equation}
From these definitions, we can immediately verify:
\begin{enumerate}
    \item $\textup{dim}[S] \geq j$ because $S$ has the same dimensionality as $J$. 
    \item For any $\mathbf{x} \in S$ and $\mathbf{y} \in J^{\prime}$ we must have
    \begin{equation}
        \label{eq:second_fact}
        \left<\mathbf{x}, \mathbf{y}\right> 
        = \left<\mathbf{P}_{\mathbf{U}} \mathbf{x}, \mathbf{y}\right> = \left<\mathbf{x}, \mathbf{P}_{\mathbf{U}} \mathbf{y}\right> = 0
    \end{equation}
    because $\mathbf{x}$ lies in the range of $\mathbf{U}$ and is orthogonal to $\mathbf{P}_{\mathbf{U}} J^{\prime}$.
    \item For any $\mathbf{x} \in H_k$ and
    $\mathbf{y} \in \textup{range}[\mathbf{U}] \cap Q_k^{\bot}$ we must have
    \begin{equation}
    \label{eq:third_fact}
    \left<\mathbf{x}, \mathbf{y}\right> 
    = \left<\mathbf{x}, \mathbf{P}_{\mathbf{U}} \mathbf{y}\right> = \left<\mathbf{P}_{\mathbf{U}} \mathbf{x}, \mathbf{y}\right> = 0,
    \end{equation}
    because
    $\mathbf{y}$ lies in the range of $\mathbf{U}$ and is orthogonal to $Q_k$.
\end{enumerate}

We next bound the distance between $S$ and $J$ by calculating
\begin{align}
    \label{eq:combo}
    \min_{\mathbf{x} \in S} 
    \frac{\mathbf{x}^{\ast} \mathbf{P}_{J} \mathbf{x}}{\left\lVert \mathbf{x} \right\rVert^2}
    & = \min_{\mathbf{x} \in S} 
    \frac{\mathbf{x}^{\ast} (\mathbf{P}_{J} + \mathbf{P}_{J^{\prime}}) \mathbf{x}}{\left\lVert \mathbf{x} \right\rVert^2} \\
    \label{eq:combo2}
    & = \min_{\mathbf{x} \in S} 
    \frac{\mathbf{x}^{\ast} \mathbf{P}_{H_k} \mathbf{x}}{\left\lVert \mathbf{x} \right\rVert^2} \\
    \label{eq:combo3}
    & = \min_{\mathbf{x} \in S} 
    \frac{\left\lVert \mathbf{P}_{H_k} \mathbf{x}\right\rVert^2}{\left\lVert \mathbf{x} \right\rVert^2} \\
    \label{eq:combo4}
    &\geq \min_{\mathbf{x} \in Q_k} 
    \frac{\left\lVert \mathbf{P}_{H_k} \mathbf{x}\right\rVert^2}{\left\lVert \mathbf{x} \right\rVert^2},
\end{align}
where \eqref{eq:combo} uses \eqref{eq:second_fact},
\eqref{eq:combo2} uses the fact that 
$J$ and $J^{\prime}$ provide an orthogonal decomposition of $H_k$,
\eqref{eq:combo3} uses the fact that $\mathbf{P}_{H_k} = \mathbf{P}_{H_k}^2$,
and \eqref{eq:combo4} uses the fact that $S \subseteq Q_k$.

To further simplify \eqref{eq:combo4}, consider orthogonal matrices $\mathbf{H}, \mathbf{Q} \in \mathbb{R}^{n \times k}$ with range $H_k, Q_k$.
Using the variational characterization of the minimum squared singular value, observe
\begin{equation}
    \sigma_{\min}(\mathbf{H}^{\ast} \mathbf{Q})^2
    = \min_{\mathbf{x} \in \mathbb{R}^k} \frac{\left\lVert \mathbf{H}^{\ast} 
    \mathbf{Q}
    \mathbf{x}\right\rVert^2}{\left\lVert \mathbf{x} \right\rVert^2}
    = \min_{\mathbf{x} \in Q_k} \frac{\left\lVert \mathbf{H}^{\ast} \mathbf{x}\right\rVert^2}{\left\lVert \mathbf{x} \right\rVert^2}
    = \min_{\mathbf{x} \in Q_k} 
    \frac{\left\lVert \mathbf{P}_{H_k} \mathbf{x}\right\rVert^2}{\left\lVert \mathbf{x} \right\rVert^2}
\end{equation}
and likewise
\begin{equation}
    \sigma_{\min}(\mathbf{Q}^{\ast} \mathbf{H})^2
    = \min_{\mathbf{x} \in \mathbb{R}^k} \frac{\left\lVert \mathbf{Q}^{\ast} 
    \mathbf{H}
    \mathbf{x}\right\rVert^2}{\left\lVert \mathbf{x} \right\rVert^2}
    = \min_{\mathbf{x} \in H_k} \frac{\left\lVert \mathbf{Q}^{\ast} \mathbf{x}\right\rVert^2}{\left\lVert \mathbf{x} \right\rVert^2}
    = \min_{\mathbf{x} \in H_k} 
    \frac{\left\lVert \mathbf{P}_{Q_k} \mathbf{x}\right\rVert^2}{\left\lVert \mathbf{x} \right\rVert^2}.
\end{equation}
This enables us to calculate
\begin{align}
    \label{eq:uses1}
    \min_{\mathbf{x} \in Q_k} 
    \frac{\left\lVert \mathbf{P}_{H_k} \mathbf{x}\right\rVert^2}{\left\lVert \mathbf{x} \right\rVert^2} &= \min_{\mathbf{x} \in H_k} 
    \frac{\left\lVert \mathbf{P}_{Q_k} \mathbf{x}\right\rVert^2}{\left\lVert \mathbf{x} \right\rVert^2} \\
    \label{eq:uses2}
    &= 
    \min_{\mathbf{x} \in H_k} 
    \frac{\left\lVert (\mathbf{P}_{Q_k}
    + \mathbf{P}_{\mathbf{U}} - \mathbf{P}_{Q_k})
    \mathbf{x}\right\rVert^2}{\left\lVert \mathbf{x} \right\rVert^2} \\
    \label{eq:uses3}
    &= \min_{\mathbf{x} \in H_k} 
    \frac{\left\lVert \mathbf{P}_{\mathbf{U}} \mathbf{x}\right\rVert^2}{\left\lVert \mathbf{x} \right\rVert^2} \\
    \label{eq:uses4}
    &= \cos^2 \theta,
\end{align}
where \eqref{eq:uses1} uses the fact that $\sigma_{\min}(\mathbf{H}^{\ast} \mathbf{Q}) = \sigma_{\min}(\mathbf{Q}^{\ast} \mathbf{H})$, \eqref{eq:uses2} uses \eqref{eq:third_fact},
\eqref{eq:uses3} follows directly from \eqref{eq:uses2},
and \eqref{eq:uses4} uses the definition of $\cos^2 \theta$.

Finishing the proof, for any $\mathbf{x} \in S$, we can decompose
\begin{equation}
\label{eq:finale1}
    \mathbf{x}
    =
    \mathbf{P}_J \mathbf{x}
    + \mathbf{P}_{H_k^{\bot}} \mathbf{x},
\end{equation}
where we have used
the fact that $J$, $J^{\prime}$ and $H_k^{\bot}$ form an orthogonal decomposition of $\mathbb{R}^n$,
and $\mathbf{x}$ must be orthogonal to $J^{\prime}$.
By the above calculations,
\begin{equation}
\label{eq:finale2}
    \min_{\mathbf{x} \in S} \frac{\mathbf{x}^{\ast} \mathbf{P}_J \mathbf{x}}{\left\lVert \mathbf{x} \right\rVert^2} \geq \cos^2 \theta,
    \qquad
    \max_{\mathbf{x} \in S} \frac{\mathbf{x}^{\ast} \mathbf{P}_{H_k^{\bot}} \mathbf{x}}{\left\lVert \mathbf{x} \right\rVert^2} \leq \sin^2 \theta.
\end{equation}
Additionally, for any $i \geq 0$, we observe that $\mathbf{P}_J \mathbf{A}^i |_{J}$
has its spectrum contained in $[\lambda_j^i, \infty)$
and $\mathbf{P}_{H_k^{\bot}} \mathbf{A}^i |_{H_k^{\bot}}$
has its spectrum contained in
$[-|\lambda_{\beta(k+1)}|^i, |\lambda_{\beta(k+1)}|^i]$, so that
\begin{equation}
\label{eq:finale3}
    \min_{\mathbf{x} \in \mathbb{R}^n}
    \frac{\mathbf{x}^{\ast}
    \mathbf{P}_J \mathbf{A}^i \mathbf{P}_J \mathbf{x}}
    {\mathbf{x}^{\ast} \mathbf{P}_J \mathbf{x}} \geq \lambda_j^i,
    \qquad
    \max_{\mathbf{x} \in \mathbb{R}^n}
    \left|
    \frac{\mathbf{x}^{\ast}
    \mathbf{P}_{H_k^{\bot}} \mathbf{A}^i \mathbf{P}_{H_k^{\bot}} \mathbf{x}}{\mathbf{x}^{\ast}
    \mathbf{P}_{H_k^{\bot}}
    \mathbf{x}}
    \right|
    \leq \left|\lambda_{\beta(k+1)}\right|^i.
\end{equation}
Lastly, using \eqref{eq:finale1}, \eqref{eq:finale2}, and \eqref{eq:finale3}, we conclude
\begin{align}
    & \quad \min_{\mathbf{x} \in S} \frac{\mathbf{x}^{\ast} \mathbf{A}^{i+1} \mathbf{x}}{\mathbf{x}^{\ast} \mathbf{A}^{i} \mathbf{x}} \\
    &= \min_{\mathbf{x} \in S} \frac{\mathbf{x}^{\ast} \mathbf{P}_J \mathbf{A}^{(i+1)} \mathbf{P}_J \mathbf{x} + \mathbf{x}^{\ast} \mathbf{P}_{H_k^{\bot}} \mathbf{A}^{i+1} \mathbf{P}_{H_k^{\bot}} \mathbf{x}}
    {\mathbf{x}^{\ast} \mathbf{P}_J \mathbf{A}^{i} \mathbf{P}_J \mathbf{x}
    + \mathbf{x}^{\ast} \mathbf{P}_{H_k^{\bot}} \mathbf{A}^{i} \mathbf{P}_{H_k^{\bot}} \mathbf{x}} \\
    &\geq \frac{\lambda_j^{i+1} \cos^2 \theta - \sin^2 \theta \left|\lambda_{\beta(k+1)}\right|^{i+1}}{\lambda_j^i \cos^2 \theta + \sin^2 \theta \left|\lambda_{\beta(k+1)}\right|^i} \\
    &= \frac{1 - \tan^2 \theta R^{i+1}}
    {1 + \tan^2 \theta R^i} \lambda_{\beta(j)}.
\end{align}
\end{proof}

The error bound \eqref{eq:explicit} is a sharper version of the error bound presented by G. W. Stewart \cite{Stewart1969}, and it has the advantage of showing clearly how the error depends on the range of starting matrix $\mathbf{U}$. 
The proof presented here is also shorter and more direct than the proof by Stewart,
taking advantage of modern strategies for manipulating the min--max principle that were originally introduced by A. Knyazev in \cite{knyazev1985sharp}.

We can simplify the error bound \eqref{eq:explicit} slightly (at the cost of reducing some of the sharpness) by taking ratios and manipulating terms to yield
\begin{equation}
\label{eq:simpler}
     0 \leq \frac{\lambda_j - \mathbf{\Lambda}^{(i)}_{jj}}{\lambda_j} \leq 2 \tan^2\theta \left| \frac{\lambda_{\beta(k+1)}}{\lambda_j} \right|^i,
     \qquad 1 \leq j \leq k.
\end{equation}
which emphasizes that the relative eigenvalue error decays exponentially fast,
in proportion to the relative eigenvalue gap
$|\lambda_{\beta(k+1)} \slash \lambda_j|^i$.
As a consequence of \eqref{eq:simpler}, we anticipate that deterministic subspace iteration is most successful at estimating the first few eigenvalues and
the accuracy degrades for eigenvalues with higher indices.

Our analysis can be adapted to the standard subspace iteration with the eigenvalue estimator \eqref{eq:symmRayRitz}.
In this case, the variational characterization of the eigenvalues becomes
\begin{equation}
    \mathbf{\Lambda}_{jj}^{(i)} = \max_{\substack{S \subseteq \textup{range}[\mathbf{U}] \\ \textup{dim}[S] = j}}
    \min_{\mathbf{x} \in S} \frac{\mathbf{x}^{\ast} \mathbf{A}^{2i+1} \mathbf{x}}{\mathbf{x}^{\ast} \mathbf{A}^{2i} \mathbf{x}}, \qquad 1 \leq j \leq k, \qquad i \geq 0,
\end{equation}
which differs slightly from \eqref{eq:variational} because the matrices $\mathbf{A}^{i}$ and $\mathbf{A}^{i+1}$ are replaced by the higher powers $\mathbf{A}^{2i}$ and $\mathbf{A}^{2i+1}$.
Consequently, the rate of convergence is twice as fast, i.e.,
\begin{equation}
     0 \leq \frac{\lambda_j - \mathbf{\Lambda}^{(i)}_{jj}}{\lambda_j}
     \leq 2 \tan^2\theta \left| \frac{\lambda_{\beta(k+1)}}{\lambda_j} \right|^{2i},
     \qquad 1 \leq j \leq k.
\end{equation}
Yet despite this faster rate of convergence, we prefer not to use the standard subspace iteration with \eqref{eq:symmRayRitz} as the foundation for our randomized subspace iteration.
Indeed, Appendix \ref{sec:symmSubsp} shows that that applying iterative random sparsification to standard subspace iteration fails dramatically,
and it is better to use the non-standard subspace iteration with \eqref{eq:genEval}.

\subsection{Analysis of stochastic compression}
\label{sec:compressionAnalysis}

In randomized subspace iteration,
we enforce sparsity in the iterates by applying \emph{stochastic compression}.
Stochastic compression is a random operation that replaces a vector $\mathbf{x}$ with a random vector $\Phi(\mathbf{x})$ that is unbiased, i.e.,
$\textup{E} \Phi(\mathbf{x}) = \mathbf{x}$
and satisfies a sparsity constraint
$\left\lVert \Phi(\mathbf{x}) \right\rVert_0 \leq m$.

Because the perturbations $\mathbf{x} - \Phi(\mathbf{x})$ contribute variance to the eigenvalue estimates,
it is desirable to make these perturbations as small as possible.
In the following proposition, we show how to do this explicitly, by minimizing the mean square perturbation magnitude $\textup{E} \left\lVert \mathbf{x} - \Phi(\mathbf{x}) \right\rVert^2_2$.

\begin{proposition}{\label{prop:optimal}}
For any $\mathbf{x} \in \mathbb{R}^n$,
the solution to the minimization problem
\begin{equation}
    \min_{\mathbf{\Phi}(\mathbf{x})} 
    \left\{\textup{E} \left\lVert \mathbf{x} - \Phi(\mathbf{x}) \right\rVert^2_2 \colon \, 
    \left\lVert \Phi\left(\mathbf{x}\right) \right\rVert_0 \leq m,
    \,
    \textup{E} \Phi\left(\mathbf{x}\right) = \mathbf{x} \right\},
\end{equation}
is characterized by three properties:
\begin{enumerate}
	\item The largest-magnitude entries $|\mathbf{x}_{\alpha\left(1\right)}| \geq \cdots \geq |\mathbf{x}_{\alpha\left(k\right)}|$ are preserved exactly, i.e., $\Phi(\mathbf{x})_{\alpha\left(i\right)} = \mathbf{x}_{\alpha\left(i\right)}$.
	\item The smallest-magnitude entries $|\mathbf{x}_{\alpha\left(k+1\right)}| \geq \cdots \geq |\mathbf{x}_{\alpha\left(n\right)}|$ are randomly perturbed, i.e., $\Phi(\mathbf{x})_{\alpha\left(i\right)} = \mathbf{x}_{\alpha\left(i\right)} \slash \mathbf{p}_{\alpha\left(i\right)}$
	with probability $\mathbf{p}_{\alpha\left(i\right)}$	
	and $\Phi(\mathbf{x})_{\alpha\left(i\right)} = 0$ otherwise,
	where
	\begin{equation}
	\mathbf{p}_{\alpha\left(i\right)} = (m - k) \left|\mathbf{x}_{\alpha\left(i\right)}\right|\slash \sum_{j = k+1}^n \left|\mathbf{x}_{\alpha\left(j\right)}\right|.
	\end{equation}
	\item The number $k$ is as small as possible while ensuring $\mathbf{p}_{\alpha\left(i\right)} \leq 1$ for $i \geq k + 1$.
\end{enumerate}

\end{proposition}
\begin{proof}
For any random vector $\Phi(\mathbf{x}) \in \mathbb{R}^n$,
introduce the vector $\Psi = \Psi(\mathbf{x}) \in \left\{0, 1\right\}^n$ with entries 
\begin{equation}
\Psi_i = \mathds{1}\left\{\Phi(\mathbf{x})_i \neq 0\right\}, \qquad 1 \leq i \leq n.
\end{equation}
$\Psi_i$ is the indicator function on the event $\Phi(\mathbf{x})_i \neq 0$, that is, its value is $1$ if this event occurs and its value is $0$ otherwise.
Using the unbiasedness condition
$\mathbf{E} \Phi(\mathbf{x})_i = \mathbf{x}_i$,
compute
\begin{equation}
\label{eq:conditional}
    \textup{E}\left[\left.\Phi(\mathbf{x})_i\right| \Psi_i
    \right] = \frac{\Psi_i}{\textup{E} \Psi_i } \mathbf{x}_i, \qquad 1 \leq i \leq n.
\end{equation}
Hence, the square error $\textup{E} \left\lVert \mathbf{x} - \Phi(\mathbf{x}) \right\rVert^2_2$
can be decomposed as
\begin{align}
    & \sum_{i=1}^n \textup{E} \left| \mathbf{x}_i - \Phi(\mathbf{x})_i \right|^2 \\
    &= \sum_{i=1}^n \textup{E} \big| \mathbf{x}_i - \textup{E}[\Phi(\mathbf{x})_i| \Psi_i
    ]\big|^2 
    + \sum_{i=1}^n 
    \textup{E} \big| \textup{E} [\Phi(\mathbf{x})_i | \Psi_i ] - \Phi(\mathbf{x})_i\big|^2 \\
    \label{eq:decompose}
    &= \sum_{i=1}^n \mathbf{x}_i^2 \left(
    \frac{1}{\textup{E} \Psi_i} - 1
    \right)
    + \sum_{i=1}^n \textup{E} \left|\frac{\Psi_i}{\textup{E} \Psi_i } \mathbf{x}_i - \Phi(\mathbf{x})_i \right|^2.
\end{align}
We minimize \eqref{eq:decompose} by taking $\Phi(\mathbf{x})_i = \mathbf{x}_i \Psi_i \slash \textup{E} \Psi_i$ for $1 \leq i \leq n$,
and choosing $\textup{E} \Psi \in \left[0, 1\right]^n$ to solve
\begin{equation}
\label{eq:convex}
    \min_{\textup{E} \Psi \in \left[0, 1\right]^n} \left\{
    \sum_{i=1}^n \mathbf{x}_i^2 \left(\frac{1}{\textup{E} \Psi_i } - 1\right) \colon \, \left\lVert \textup{E} \Psi \right\rVert_1 \leq m\right\},
\end{equation}
which is a convex optimization problem with linear inequality constraints.

To solve \eqref{eq:convex}, 
let $\alpha\left(1\right), \ldots, \alpha\left(n\right)$ be a permutation that reorders the elements of $\mathbf{x}$ from largest to smallest magnitude, i.e.,
$\left|\mathbf{x}_{\alpha\left(1\right)}\right| \geq \cdots \geq \left|\mathbf{x}_{\alpha\left(n\right)}\right|$,
and introduce the Lagrangian function
\begin{equation}
    \mathcal{L}\left(\textup{E} \Psi, \eta, \boldsymbol{\lambda}\right)
    = \sum_{i=1}^n \mathbf{x}_i^2 \left(\frac{1}{\textup{E} \Psi_i } - 1\right)
    + \eta \left(\left\lVert \textup{E} \Psi \right\rVert_1 - m\right) + \sum_{i=1}^n \boldsymbol{\lambda}_i \left(\textup{E} \Psi_i - 1\right).
\end{equation}
Then, the solution to \eqref{eq:convex} must satisfy the gradient condition
$\nabla_{\textup{E} \Psi} \mathcal{L}\left(\textup{E} \Psi, \eta, \boldsymbol{\lambda}\right) = \mathbf{0}$, which leads to
\begin{equation}
    \textup{E} \Psi_i = \frac{\left|\mathbf{x}_i\right|}{\left(\eta + \boldsymbol{\lambda}_i\right)^{1 \slash 2}}, \quad 1 \leq i \leq n,
\end{equation}
as well as the complementarity condition $\boldsymbol{\lambda}_i \left(\textup{E} \Psi_i - 1\right) = 0$, which leads to
\begin{equation}
    \textup{E} \Psi_i = \frac{\left|\mathbf{x}_i\right|}{\eta^{1 \slash 2}} \quad \text{or} \quad \textup{E} \Psi_i = 1, \quad 1 \leq i \leq n.
\end{equation}
Examining the objective function \eqref{eq:convex} shows that it is best to set $\textup{E} \Psi_i = 1$ only for the \emph{largest-magnitude} entries $\mathbf{x}_i$. In summary, the largest-magnitude entries satisfy
\begin{equation}
    \textup{E} \Psi_{\alpha\left(i\right)} = 1, \quad 1 \leq i \leq k,
\end{equation}
and the smallest-magnitude entries satisfy
\begin{equation}
\label{eq:constraint}
    \textup{E} \Psi_{\alpha\left(i\right)} = \frac{\left(m - k\right) \left|\mathbf{x}_{\alpha\left(i\right)}\right|}{\sum_{j=k+1}^n \left|\mathbf{x}_{\alpha\left(j\right)}\right|}, \quad k+1 \leq i \leq n,
\end{equation}
for some parameter $1 \leq k \leq n$.
The value of the objective function \eqref{eq:convex} becomes
\begin{equation}
    \label{eq:largest}
    \sum_{i=k+1}^n \mathbf{x}_{\alpha\left(i\right)}^2 \left(\frac{\sum_{j=k+1}^n \left|\mathbf{x}_{\alpha\left(j\right)}\right|}{\left(m - k\right) \left|\mathbf{x}_{\alpha\left(i\right)}\right|} - 1\right),
\end{equation}
and a direct computation reveals \eqref{eq:largest} to be nondecreasing in $k$, whence $k$ should be taken as small as possible while ensuring that $\textup{E} \Psi_i \leq 1$ for all $1 \leq i \leq n$.
\end{proof}

Proposition \ref{prop:optimal} identifies necessary and sufficient conditions
to minimize the mean square perturbation magnitude $\textup{E} \left\lVert \mathbf{x} - \Phi(\mathbf{x}) \right\rVert^2_2$ over all compression schemes.
Under these conditions, pivotal compression is an optimal scheme,
but the multinomial compression scheme discussed in \cite{Lim2017} and \cite{Greene2019} is not optimal.
Building on this comparison, Appendix \ref{sec:pivVsys} empirically compares pivotal compression and two other schemes and finds that pivotal compression leads to the highest accuracy eigenvalue estimates.

\subsection{Analysis of randomized subspace iteration}
\label{sec:randSubspAnalysis}

For the last part of the analysis, we observe that the iterates $\mathbf{X}^{(0)}, \mathbf{X}^{(1)}, \mathbf{X}^{(2)}, \ldots$ from randomized subspace iteration are embedded in a Markov chain.
Therefore, we can apply general results governing the asymptotic behavior of Markov chains,
specifically the strong law of large numbers and Markov chain central limit theorem, which hold under a geometric ergodicity assumption \cite{meyn2012markov,jones2004markov}.
This asymptotic theory is general to Markov chains, so it does not lead to any specific \emph{a priori} bounds concerning the accuracy of randomized subspace iteration.
However, this theory does lead to qualitative predictions:
first that the estimated eigenvalues $\mathbf{\Lambda}^{(i_{\max})}$ converge to a deterministic limit as $i_{\max} \rightarrow \infty$
and second that the variation around this limit is asymptotically Gaussian with a variance that can be estimated \emph{a posteriori} from data.

We emphasize that the geometric ergodicity assumption for the Markov chain is significant.
Geometric ergodicity implies that the Markov chain `resets' itself at regular intervals,
and this condition can be violated for example if the matrix $\mathbf{A}$ is block diagonal.
As for any Markov chain, irreducibility is an important prerequisite for ergodicity~\cite{Levin2017mixing} and rules out the block diagonal case.
However, a full identification of the conditions guaranteeing geometric ergodicity is beyond the scope of the current analysis.

Recalling the definition $\mathbf{Y}^{(i)} = \mathbf{A} \Phi(\mathbf{X}^{(i)})$,
the Markov chain is defined by
\begin{multline}
    \mathbf{\Gamma}^{(i)} = \big[\mathbf{X}^{(\Delta i)}\
    \mathbf{Y}^{(\Delta i)}\
    \mathbf{X}^{(\Delta i + 1)}\
    \mathbf{Y}^{(\Delta i + 1)}\ \cdots \\ \mathbf{X}^{(\Delta (i + 1) - 1)}\
    \mathbf{Y}^{(\Delta (i + 1) - 1)}\
    \mathbf{N}^{(\Delta (i + 1) - 1)}\big],
\end{multline}
where $\Delta$ is the orthogonalization interval ($\Delta = 1000$ in our experiments).
Because the sequence $\mathbf{\Gamma}^{(0)}, \mathbf{\Gamma}^{(1)}, \ldots$ satisfies the conditional independence property
\begin{equation}
\label{eq:weak}
    \textup{Law}\left(\mathbf{\Gamma}^{(i)} |
    \mathbf{\Gamma}^{(i - 1)}, \ldots, \mathbf{\Gamma}^{(0)}\right) =
    \textup{Law}\left(\mathbf{\Gamma}^{(i)} |
    \mathbf{\Gamma}^{(i-1)}\right),
\end{equation}
it is indeed Markovian.

\begin{proposition}
\label{prop:asymptotics}
Assume the Markov chain $\mathbf{\Gamma}^{(0)}, \mathbf{\Gamma}^{(1)}, \ldots$ is geometrically ergodic with respect to a distribution $\mu$,
introduce the matrices
\begin{equation}
    \mathbf{E}^{(\infty)} = \textup{E}_{\mu}\Biggl[\frac{1}{\Delta} \sum_{i=0}^{\Delta - 1} \mathbf{U}^{\ast} \mathbf{X}^{(i)}\Biggr] \qquad \text{and} \qquad
    \mathbf{F}^{(\infty)} =
    \textup{E}_{\mu}\Biggl[\frac{1}{\Delta} \sum_{i=0}^{\Delta - 1}
    \mathbf{U}^{\ast} \mathbf{Y}^{(i)}\Biggr],
\end{equation}
and assume
$\textup{E}_{\mu}\lVert\frac{1}{\Delta} \sum_{i=0}^{\Delta - 1} \mathbf{U}^{\ast} \mathbf{X}^{(i)}\rVert_{\textup{F}}^{2 + \epsilon} < \infty$
and 
$\textup{E}_{\mu}\lVert\frac{1}{\Delta} \sum_{i=0}^{\Delta - 1}
\mathbf{U}^{\ast} \mathbf{Y}^{(i)}\rVert_{\textup{F}}^{2 + \epsilon} < \infty$ for some $\epsilon > 0$.
Lastly, assume $\mathbf{E}^{(\infty)}$ is nonsingular.
Then the estimated eigenvalues from randomized subspace iteration converge with probability one
\begin{equation}
\label{eq:convergence}
    \lim_{i_{\max} \rightarrow \infty} \mathbf{\Lambda}^{(i_{\max})} = \mathbf{\Lambda}^{(\infty)},
\end{equation}
where $\mathbf{\Lambda}^{(\infty)}$ is the
the solution to the generalized eigenvalue problem
\begin{equation}
    \mathbf{F}^{(\infty)} \mathbf{W}^{(\infty)} \mathbf{\Lambda}^{(\infty)}
    = \mathbf{E}^{(\infty)} \mathbf{W}^{(\infty)}.
\end{equation}

Additionally, suppose $\mathbf{\Lambda}_{jj}^{(\infty)}$ is a simple eigenvalue with corresponding left eigenvector $\mathbf{z}^{(\infty)}_j$ and right eigenvector $\mathbf{w}^{(\infty)}_j$, set
\begin{equation}
\label{eq:var_expr}
    f_j^{(i)} = \mathbf{z}_j^{(\infty)} \left( \mathbf{U}^{\ast} \mathbf{X}^{(i)} - \mathbf{\Lambda}_{jj}^{(\infty)} \mathbf{U}^{\ast} \mathbf{Y}^{(i)} \right) \mathbf{w}_j^{(\infty)},
\end{equation}
and define the asymptotic variance
\begin{equation}
\label{eq:as_defined}
    \sigma^2_j = \lim_{i_{\max} \rightarrow \infty} 
    \frac{1}{i_{\max}} \sum_{i,k=0}^{i_{\max} - 1} 
    \textup{Cov}_{\mu}\left[f_j^{(i)}, f_j^{(k)}\right].
\end{equation}
Then, the deviations $\mathbf{\Lambda}_{jj}^{(i_{\max})} - \mathbf{\Lambda}_{jj}^{(\infty)}$ are asymptotically Gaussian
\begin{equation}
\label{eq:clt}
    \sqrt{i_{\max}}\left(\mathbf{\Lambda}_{jj}^{(i_{\max})} - \mathbf{\Lambda}_{jj}^{(\infty)}\right) \stackrel{\mathcal{D}}{\rightarrow} \mathcal{N}\left(0, \sigma^2_j\right)
\end{equation}
and satisfy the law of the iterated logarithm with probability one
\begin{multline}
\label{eq:lil}
    -\sqrt{2} \sigma_j = \liminf_{i_{\max} \rightarrow \infty} \sqrt{\frac{i_{\max}}{\log \log i_{\max}}} \left(\mathbf{\Lambda}_{jj}^{(i_{\max})} - \mathbf{\Lambda}_{jj}^{(\infty)}\right) \\
    \leq \limsup_{i_{\max} \rightarrow \infty} \sqrt{\frac{i_{\max}}{\log \log i_{\max}}} \left(\mathbf{\Lambda}_{jj}^{(i_{\max})} - \mathbf{\Lambda}_{jj}^{(\infty)}\right)
    = \sqrt{2} \sigma_j.
\end{multline}
\end{proposition}
\begin{proof}
Using the geometric ergodicity assumption and strong law of large numbers for Markov chains \cite{meyn2012markov}, we confirm that almost surely
\begin{align}
\label{eq:mat_one}
    & \mathbf{E}^{(i_{\max})} =
    \frac{1}{i_{\max} - i_{\min}}
    \sum_{i=i_{\min}}^{i_{\max} - 1}
    \mathbf{U}^{\ast} \mathbf{X}^{(i)}
    \rightarrow \mathbf{E}^{(\infty)}, \\
    \label{eq:mat_two}
    & \mathbf{F}^{(i_{\max})} =
    \frac{1}{i_{\max} - i_{\min}}
    \sum_{i=i_{\min}}^{i_{\max} - 1}
    \mathbf{U}^{\ast} \mathbf{Y}^{(i)}
    \rightarrow
    \mathbf{F}^{(\infty)},
\end{align}
as $i_{\max} \rightarrow \infty$.
Moreover, \eqref{eq:mat_one}, \eqref{eq:mat_two}, and the continuity of eigenvalues in the nonsingular generalized eigenvalue problem \cite{Stewart2001vol2} together imply that $\mathbf{\Lambda}^{(i_{\max})} \rightarrow \mathbf{\Lambda}^{(i_{\infty})}$ almost surely, confirming \eqref{eq:convergence}.

To prove the asymptotic Gaussianity
and law of the iterated logarithm, we use the first-order perturbation theory for simple eigenvalues (see e.g., \cite{Stewart2001vol2,greenbaum2020first}), which reveals that
\begin{multline}
\label{eq:linear}
    \mathbf{\Lambda}_{jj}^{(i)} - \mathbf{\Lambda}_{jj}^{(\infty)} \\
    = \mathbf{z}_j^{(\infty)} \left(\mathbf{E}^{(i)} - \mathbf{\Lambda}_{jj}^{(\infty)} \mathbf{F}^{(i)}\right) \mathbf{w}_j^{(\infty)}
    + \mathcal{O}\left(\mathbf{E}^{(i)} - \mathbf{E}^{(\infty)}\right)^2
    + \mathcal{O}\left(\mathbf{F}^{(i)} - \mathbf{F}^{(\infty)}\right)^2
\end{multline}
as $i_{\max} \rightarrow \infty$.
By the Markov chain law of the iterated logarithm (see \cite{oodaira1971law,jones2004markov}),
the quantities $\mathbf{E}^{(i)} - \mathbf{E}^{(\infty)}$
and $\mathbf{F}^{(i)} - \mathbf{F}^{(\infty)}$
are almost surely $\mathcal{O}(\sqrt{i_{\max} \slash \log \log i_{\max}})$ as $i_{\max} \rightarrow \infty$.
Moreover, by the Markov chain central limit theorem and law of the iterated logarithm (see \cite{ibragimov1962some,oodaira1971law,jones2004markov}), we find that
\begin{equation}
\label{eq:clt_basic}
    \sqrt{i_{\max}} \mathbf{z}_j^{(\infty)} \left(\mathbf{E}^{(i)} - \mathbf{\Lambda}_{jj}^{(\infty)} \mathbf{F}^{(i)}\right) \mathbf{w}_j^{(\infty)}
    \stackrel{\mathcal{D}}{\rightarrow} \mathcal{N}(0, \sigma^2_j),
\end{equation}
and almost surely
\begin{multline}
\label{eq:lil_basic}
    -\sqrt{2} \sigma_j = \liminf_{i_{\max} \rightarrow \infty} \sqrt{\frac{i_{\max}}{\log \log i_{\max}}} \mathbf{z}_j^{(\infty)} \left(\mathbf{E}^{(i)} - \mathbf{\Lambda}_{jj}^{(\infty)} \mathbf{F}^{(i)}\right) \mathbf{w}_j^{(\infty)} \\
    \leq \limsup_{i_{\max} \rightarrow \infty} \sqrt{\frac{i_{\max}}{\log \log i_{\max}}} \mathbf{z}_j^{(\infty)} \left(\mathbf{E}^{(i)} - \mathbf{\Lambda}_{jj}^{(\infty)} \mathbf{F}^{(i)}\right) \mathbf{w}_j^{(\infty)}
    = \sqrt{2} \sigma_j,
\end{multline}
where the asymptotic variance $\sigma^2_j$ is defined in \eqref{eq:as_defined}.
Using Slutsky's lemma \cite{gut2013probability}, we confirm \eqref{eq:clt} and \eqref{eq:lil}.
\end{proof}

As a consequence of Proposition \ref{prop:asymptotics},
the quality of each eigenvalue estimate $\mathbf{\Lambda}_{jj}^{(i_\text{max})}$ is mainly determined by the \emph{asymptotic variance} and the \emph{asymptotic bias}.
The asymptotic variance is the stochastic variability of $\mathbf{\Lambda}_{jj}^{(i_\text{max})}$ in the limit as $i_{\max} \rightarrow \infty$,
as defined by
\begin{equation}
    \sigma_j^2 = \lim_{i_{\max} \rightarrow \infty} i_{\max} \textup{Var}[\mathbf{\Lambda}_{jj}^{(i_{\max})}].
\end{equation}
The asymptotic bias is the difference between the estimate after infinitely many iterations and the exact eigenvalue,
as defined by
\begin{equation}
    b_j = \bigl|\lim\nolimits_{i_\text{max} \rightarrow \infty} \mathbf{\Lambda}_{jj}^{(i_\text{max})} - \lambda_j\bigr|.
\end{equation}
Empirically, we find that increasing $m$ in compression operations reduces the standard error $\sigma_j^2$ and the asymptotic bias $b_j$. 
In contrast, averaging over more iterations (i.e.~increasing ${i_\text{max}}$) reduces only the variance.

Lastly, by taking advantage of the explicit expression \eqref{eq:var_expr},
we can compute the standard error in each eigenvalue estimate as follows.
We let ${\mathbf{w}}_j$ and $\mathbf{z}_j$
denote the left and right generalized eigenvectors corresponding to the eigenvalue $\mathbf{\Lambda}_{jj}^{(i_\text{max})}$ in \eqref{eq:linEst}
and define the scalar-valued time series
\begin{equation}
f^{(i)}_j = {\mathbf{z}}_j^* \left(\mathbf{U}^* \mathbf{Y}^{(i)} - \mathbf{\Lambda}_{jj}^{(i_\text{max})} \mathbf{U}^* \mathbf{X}^{(i)} \right) {\mathbf{w}}_j.
\end{equation}
Then, from \eqref{eq:as_defined} the variance in $\mathbf{\Lambda}_{jj}^{(i_\text{max})}$ is well-approximated by
\begin{equation}
\label{eq:varEstimate}
\textup{Var}[\mathbf{\Lambda}_{jj}^{(i_\text{max})}]
\approx \frac{1}{(i_\text{max} - i_\text{min})^2}
\sum_{\substack{i_{\min} \leq i,k \leq i_{\max} - 1 \\
\left|i - k\right| \leq \tau
}}
f_j^{(i)} f_j^{(k)},
\end{equation}
where the truncation threshold $\tau$ is chosen large enough
that the correlations involving $f_j^{(i)}$ and $f_j^{(k)}$ are negligibly small for any $|i - k| > \tau$ \cite{Sokal1997}.
These \emph{a posteriori} variance formulas enable estimating the standard error terms $\sigma_j^{(i_{\max})}$,
which can be useful for measuring and reducing stochastic errors in eigenvalue estimates.

\section{Numerical experiments}
\label{sec:results}
We next assess the performance of our method by applying it to the full configuration interaction (FCI) problem from quantum chemistry.
Quantum chemistry methods like FCI can be used to predict the properties of molecules (e.g. their atomic geometries, behavior in chemical reactions, or response to stimulation by light).
In most cases, estimating eigenvalues to within $1$ m$E_\text{h}$ is necessary for chemical accuracy~\cite{Bogojeski2020}, where m$E_\text{h}$ denotes the milliHartree unit of energy.

The FCI Hamiltonian matrix $\mathbf{H}$ encodes the physics of interacting electrons in a field of fixed nuclei.
It is expressed in a basis of Slater determinants, each representing a configuration of $N$ electrons in $M$ single-electron orbitals, here taken to be canonical Hartree--Fock orbitals.
As a normalization, we subtract the Hartree--Fock energy from all diagonal elements.
As $N$ and $M$ increase, the dimensionality of $\mathbf{H}$ scales as $\mathcal{O}\bigl({M \choose N \slash 2}^2\bigr)$.
The elements of $\mathbf{H}$ are given by the Slater--Condon rules \cite{Slater1929, Condon1930}, which also determine its sparsity structure; only $\mathcal{O}(N^2 M^2)$ elements per column are nonzero.

For many molecules, symmetry can be used to construct a basis in which the Hamiltonian matrix is block diagonal, effectively reducing the dimensionality.
Here, we consider only the spatial (i.e.~point-group) symmetry of the nuclei, although other symmetries (e.g. electron spin or angular momentum) could also be used.
Many randomized methods leverage this block-diagonal structure to estimate multiple eigenvalues by calculating the dominant eigenvalue for each block independently.
Here we focus on the more challenging task of estimating multiple eigenvalues within the single block containing the lowest-energy eigenvalue.

Table \ref{tab:params} lists the parameters defining the FCI matrix for each molecular system considered here.
We used the PySCF software \cite{Sun2018} to calculate matrix elements, point-group symmetry labels, and reference eigenvalues to which results are compared.
The dimensions of these problems are small enough that their eigenvalues can be obtained by standard deterministic iterative methods.
For such small problems, deterministic methods are often more efficient than the randomized approaches presented here, especially considering the slow convergence of deterministic subspace iteration relative to better-performing schemes like Jacobi--Davidson~\cite{Davidson1975, Sleijpen1996}.
Nevertheless, we focus here on these small problems in order to demonstrate the convergence behavior of our randomized approach with comparisons to exact eigenvalues.
In a companion work, we demonstrate that algorithms based on our randomization approach can yield highly accurate estimates of the excited-state energies~\cite{greene2022full} in high-dimensional FCI problems that cannot be solved using conventional methods.

The lowest-energy eigenvalues of $\mathbf{H}$ are of greatest chemical interest, so we apply our algorithm to $\mathbf{A} = \mathbf{I} - \varepsilon \mathbf{H}$ instead of $\mathbf{H}$ itself.
Eigenvalues $E_j$ of $\mathbf{H}$ are related to those $\lambda_j$ of $\mathbf{A}$ as $E_j = \varepsilon^{-1}(1 - \lambda_j)$.
For small enough $\varepsilon > 0$,
the lowest-energy eigenvalues of $\mathbf{H}$ correspond to the dominant eigenvalues of $\mathbf{A}$.
We use $\varepsilon = 10^{-6}~\text{m}E_\text{h}^{-1}$ for all systems.

\begin{table}
\caption{Parameters used in numerical calculations of FCI eigenvalues. $N$ and $M$ denote the number of active electrons and orbitals, respectively, for each system, as determined by the choice of single-electron basis. $n_\text{FCI}$ is the dimension of the block of the FCI matrix containing the lowest-energy eigenvalue.}
\begin{center}
\begin{tabular}{c | c | c | c | c}
& Nuclear & Single-electron &  &  \\
System & separation &  basis & $(N, M)$ & $n_\text{FCI} / 10^6$ \\ \hline
Ne & - & aug-cc-pVDZ & $(8, 22)$ & 6.69  \\
equilibrium \ce{C2} & $1.27273$ \AA & cc-pVDZ & $(8, 26)$ & 27.9\\
stretched \ce{C2} & $2.22254$ \AA & cc-pVDZ & $(8, 26)$ & 27.9\\
\end{tabular}
\label{tab:params}
\end{center}
\end{table}

\begin{figure}
    \centering
    \includegraphics[width=0.8\linewidth]{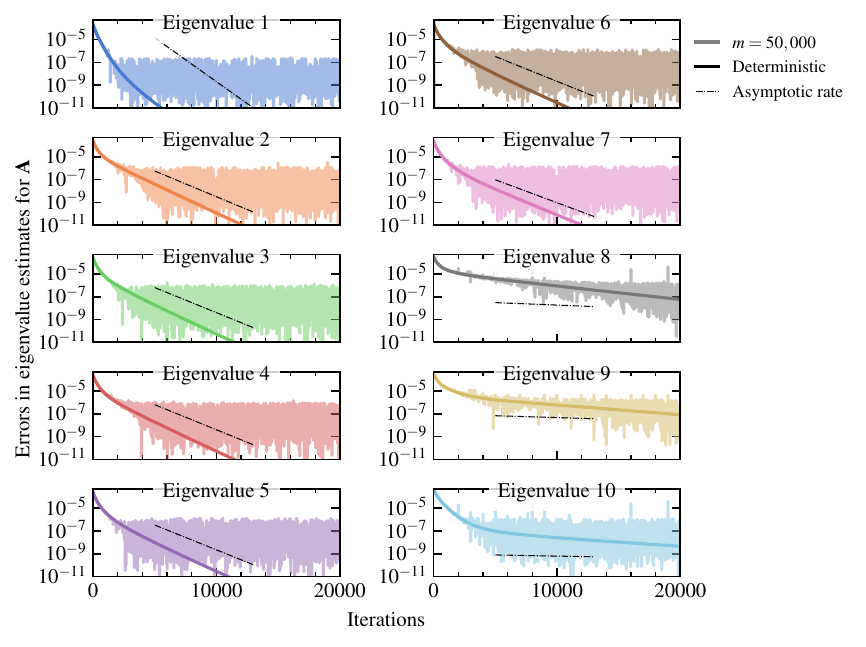}
    \caption{Errors in instantaneous Ritz values (relative to exact) for the matrix $\mathbf{A} = \mathbf{I} - \varepsilon \mathbf{H}$ corresponding to the Ne atom, obtained by applying deterministic subspace iteration and our randomized subspace iteration with $m=50,000$ nonzero elements per column. Dash-dotted lines indicate asymptotic convergence rates proportional to the relative eigenvalue gap $|\lambda_{(k + 1)} / \lambda_j|^i$ (\ref{eq:simpler}). Eigenvalues are ordered such that Eigenvalue 1 is the dominant eigenvalue of $\mathbf{A}$.}
    \label{fig:deterministic}
\end{figure}

\subsection{Randomized results for the Ne atom}
\label{sec:detNe}

In order to demonstrate the convergence behavior of randomized subspace iteration (Algorithm \ref{alg:randSubsp}), we apply it to the matrix $\mathbf{A} = \mathbf{I} - \varepsilon \mathbf{H}$ for the Ne atom in the aug-cc-pVDZ basis.
The dimensions of this matrix are approximately 7 million by 7 million,
and we set the sparsity parameter to be much smaller, specifically $m =$ 50,000.
We construct the matrix $\mathbf{U}$ of initial eigenvector estimates by diagonalizing $\mathbf{H}$ in a subspace of restricted orbital occupations (i.e.~a restriction from 22 to just 10 orbitals), leading to a sparse representation with approximately 5000 nonzero elements per column.
For comparison purposes we also apply the deterministic subspace iteration from Section \ref{sec:detSubsp}.

For our initial comparison between the randomized and deterministic schemes, we do not apply the averaging techniques from Section \ref{sec:randSubsp}. Instead, Figure \ref{fig:deterministic} presents the errors in the instantaneous Ritz values obtained by solving \eqref{eq:genEval} after each iteration.
The dash-dotted lines indicate the asymptotic convergence rates for the deterministic eigenvalue estimates, as derived in Section \ref{sec:deterministicAnalysis}.
During an initialization period of $\sim 2 \times 10^3$ iterations, the instantaneous Ritz values for the randomized and deterministic schemes are similar and converge \emph{more quickly} than the asymptotic rates.
However, after the initialization period,
the deterministic eigenvalue estimates converge more slowly, and the randomized eigenvalue estimates settle into a stationary distribution, in which they exhibit high-frequency random fluctuations due to the compression operations.

\begin{figure}
    \centering
    \includegraphics[scale=0.75]{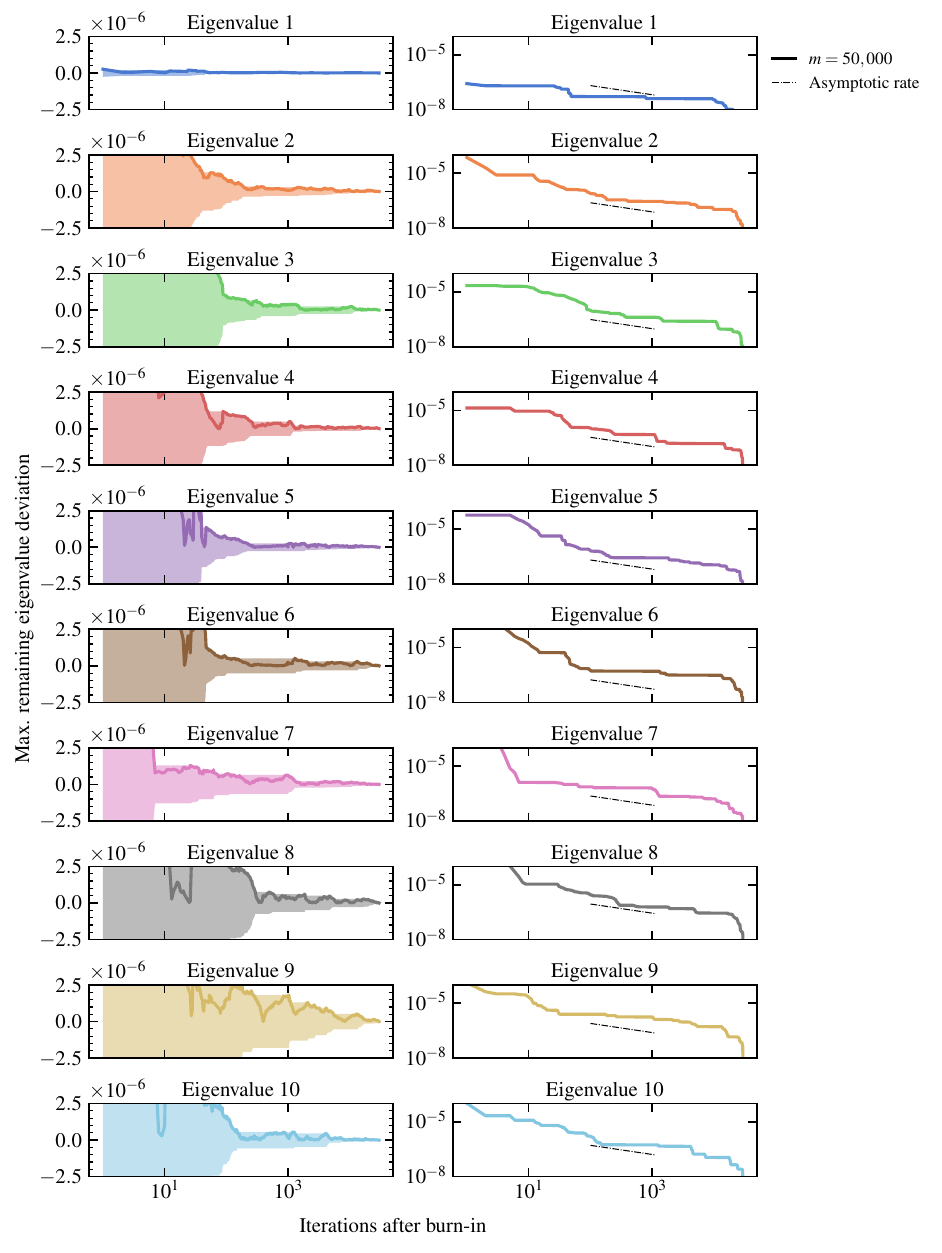}
    \caption{Convergence of the averaged eigenvalue estimates for the matrix $\mathbf{A} = \mathbf{I} - \varepsilon \mathbf{H}$ corresponding to the Ne atom, obtained with $m=50,000$ with respect to the iterations after the burn-in time, $i_\text{min}$.}
    \label{fig:meanConv}
\end{figure}

Next, we apply the averaging technique \eqref{eq:linEst} to improve the eigenvalue estimates for the randomized scheme.
Figure \ref{fig:meanConv} shows what happens when we start averaging after a burn-in time of $i_{\min} = 2 \times 10^4$ iterations.
As the number of iterations increases, the averaged eigenvalue estimates fluctuate less and less, ultimately converging to an asymptotic limit.
The widths of the shaded regions indicate the maximum future fluctuations around the asymptotic limit, that is,
\begin{equation}
    \textup{width}_j^i = \max_{i^{\prime} \geq i} |\mathbf{\Lambda}^{(i^{\prime})}_{jj} - \mathbf{\Lambda}_{jj}^{(\infty)}|
\end{equation}
Lastly, the right panel of Figure \ref{fig:meanConv} shows the maximum future fluctuations on a log-log scale to reveal further details.
Consistent with the expected convergence for averages of a Markov chain (see Proposition \ref{prop:asymptotics}), the maximum future fluctuations decrease at a rate that is asymptotically $1 \slash \sqrt{i}$ (more precisely $\sqrt{\log \log i \slash i}$) as $i \rightarrow \infty$,
as indicated by the dash-dotted line.

\begin{figure}
    \centering
    \includegraphics[scale=0.75]{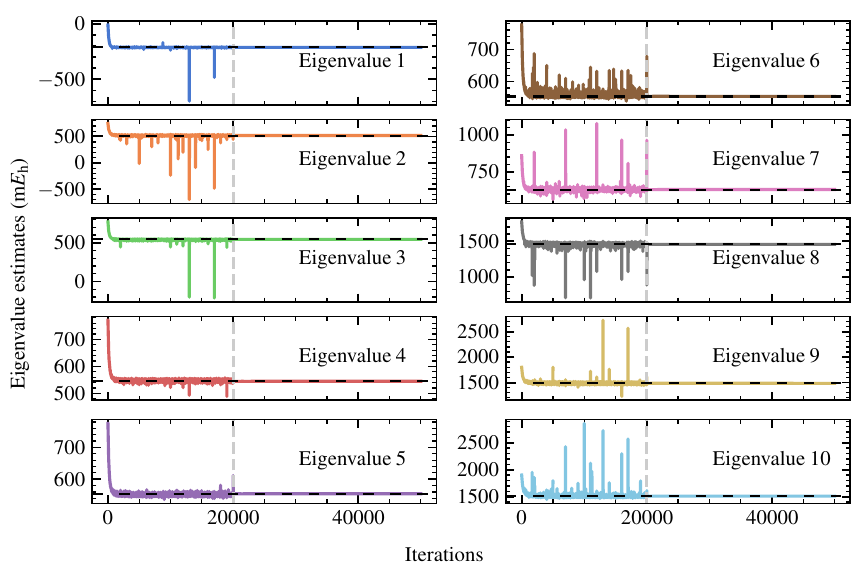}
    \caption{Eigenvalue estimates from our randomized subspace iteration. After an initial equilibration period (20,000 iterations), estimates of the first ten eigenvalues of $\mathbf{H}$ for the Ne atom converge to exact eigenvalues (horizontal dashed lines) to within small asymptotic biases when iterate matrices are compressed to $m=10,000$ nonzero elements per column.}
    \label{fig:NeEst}
\end{figure}



Fig. \ref{fig:NeEst} presents eigenvalue estimates for this same Ne system, obtained by performing stochastic compressions with only $m=10,000$ nonzero elements per column in each iteration.
We show the instantaneous Ritz values for the first 20,000 iterations,
after which we begin averaging (i.e.~with $i_\text{min} = 20,000$). Eigenvalues of the averaged matrices show substantially decreased fluctuations.

\begin{figure}
    \centering
    \includegraphics[scale=0.75]{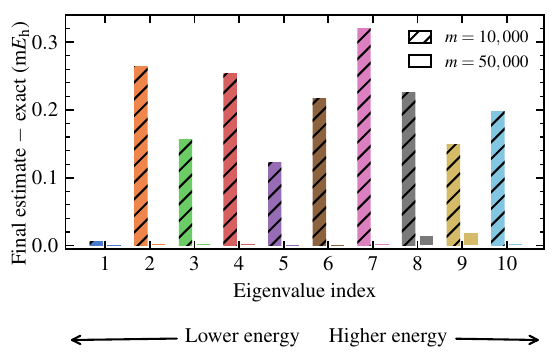}
    \caption{Errors in final eigenvalue estimates for $\mathbf{H}$ from our randomized subspace iteration after 50,000 iterations. Eigenvalue estimates obtained with $m=50,000$ exhibit less error than those with $m=10,000$. All estimates agree with exact eigenvalues to within 0.35~m$E_\text{h}$.}
    \label{fig:NeErr}
\end{figure}

Figure \ref{fig:NeErr} shows differences between estimated and exact eigenvalues after 50,000 iterations.
With $m=10,000$, these differences are less than 0.32 m$E_\text{h}$ and can be attributed to the asymptotic biases in our method.
All standard errors are less than $0.004$ m$E_\text{h}$.
Retaining more nonzero elements in compression operations $(m=50,000)$ reduces the asymptotic biases to $<0.02$  m$E_\text{h}$ and the standard errors to $<0.0003$ m$E_\text{h}$.


\begin{table}
\caption{Exact eigenvalues (in m$E_\text{h}$) for the equilibrium and stretched \ce{C2} systems, and those calculated by our randomized subspace iteration with $m=150,000$ or $m=200,000$, respectively.}
\begin{center}
\begin{tabular}{r r c r r}
\multicolumn{2}{c}{equilibrium \ce{C2} ($m=150,000$)} & & \multicolumn{2}{c}{stretched \ce{C2} ($m=200,000$)} \\
Exact  &  Estimated & & Exact  &  Estimated \\ \hline
$-343.40$ &  $-343.39$ & \hspace{1em} & $-346.33$ & $-346.32$ \\
$-264.67$ & $ -264.67$ & \hspace{1em} & $-343.50$ & $-343.50$ \\
$-254.34$ & $-254.34$ & \hspace{1em} & $-322.30$ & $-322.33$ \\
$-155.40$ & $-155.37$ & \hspace{1em} & $-298.45$ & $-298.45$ \\ 
$-104.86$ & $-104.84$ & \hspace{1em} & $-251.24$ & $-251.22$ \\
$-85.17$ & $-85.15$ & \hspace{1em} & $-250.07$ & $-250.05$ \\
$-48.91$ & $-48.92$ & \hspace{1em} & $-241.71$  & $-241.72$ \\
$-19.33$ & $-19.35$ & \hspace{1em} & $-208.54$ & $-208.15$ \\
$13.64$ & $13.63$ & \hspace{1em} & $-201.02$ & $-199.05$ \\
$26.44$ & $26.53$ & \hspace{1em} & $-190.49$ & $-190.01$ \\
\end{tabular}
\label{tab:C2data}
\end{center}
\end{table}


\subsection{Randomized results for the carbon dimer}
\label{sec:C2results}

In order to evaluate the performance of our randomized subspace iteration for higher-dimensional problems, we applied it to the carbon dimer molecule \ce{C2}, with computational details given in Table \ref{tab:params} and results given in Table \ref{tab:C2data}.
For both equilibrium and stretched \ce{C2},
we constructed the matrix $\mathbf{U}$ of initial eigenvector estimates by diagonalizing $\mathbf{H}$ in a subspace of restricted orbital occupations (i.e.~a restriction from 26 to just 9 orbitals).

Randomized subspace iteration performs well for \ce{C2} at its equilibrium bond length. 
With iterates compressed to $m=150,000$ nonzero elements per column, estimates of the ten lowest-energy eigenvalues differ from exact eigenvalues by less than 0.09 m$E_\text{h}$.
Standard errors after 50,000 iterations (choosing $i_\text{min} = 20,000$) are less than $4 \times 10^{-5}$ m$E_\text{h}$.

The stretched \ce{C2} molecule exhibits stronger electron correlation than the equilibrium \ce{C2} molecule and therefore serves as a more rigorous test for numerical eigenvalue methods in general.
Deterministic subspace iteration converges more slowly for stretched \ce{C2} because differences between eigenvalues are smaller (Table \ref{tab:C2data}).
Our randomized subspace iteration performs worse for stretched \ce{C2} than for equilibrium \ce{C2} when applied with the same parameters ($m=150,000$).
As the iteration proceeds, random fluctuations in elements of the matrices $\mathbf{U}^* \mathbf{A} \Phi(\mathbf{X}^{(i)})$ and $\mathbf{U}^* \mathbf{X}^{(i)}$ increase, and ultimately the matrices $\mathbf{U}^* \mathbf{X}^{(i)}$ become singular.
This makes it impossible to obtain accurate eigenvalue estimates.
In contrast, these fluctuations are significantly reduced when iterates are instead compressed to $m=200,000$ nonzero elements per column.
In this case, eigenvalue estimates differ from exact eigenvalues by less than 2 m$E_\text{h}$, with standard errors less than $9 \times 10^{-5}$ m$E_\text{h}$.
A longer equilibration period ($i_\text{min} = 35,000$) is needed to obtain accurate estimates for this problem.

These results for equilibrium and stretched \ce{C2} suggest that each problem may require a minimum value of $m$ to achieve reliable convergence.
Others \cite{Spencer2012, Vigor2016} have observed similarly abrupt changes in variance with the amount of sampling in randomized methods for single eigenvalues, but in those cases the variance can in principle be reduced just by averaging over more iterations \cite{Greene2020}.
In contrast, insufficient sampling in our randomized subspace iteration precludes the estimation of $k$ eigenvalues due to numerical issues encountered when solving \eqref{eq:linEst}, an issue that cannot be remedied by including more iterations.

\section{Discussion}
\label{sec:conclusions}
Incorporating random sparsification techniques into iterative linear algebra methods can enable substantial gains in computational efficiency.
Here we present a general technique for extending existing iterative approaches with repeated random sparsification from estimating one dominant eigenvalue to multiple eigenvalues.
We evaluate its performance in the context of one specific randomization scheme,  pivotal compression.
We perform numerical experiments on FCI problems from quantum chemistry involving matrices as large as 28 million by 28 million.
Even when the number of elements retained in each iteration is less than 1\% of the dimension of the matrix, we obtain accurate eigenvalue estimates for three different systems.

Among previous randomized iterative approaches to the multiple eigenvalue problem, ours is perhaps most closely related to ``replica'' schemes introduced for quantum applications. Replica schemes use two independent sequences of randomly generated vectors to build subspaces within which the target matrix is subsequently diagonalized \cite{Overy2014, Blunt2015a, Blunt2015,Blunt2018nonlinear}.
In comparison, our method avoids the high-variance inner products of sparse, random vectors that can hinder replica approaches \cite{Blunt2015a} and results in a stable stochastic iteration that can be averaged to further reduce statistical error.

A variety of possible improvements could enable application of our method to even larger problems.
Several of these have been developed and tested in the context of the ground-state (i.e.~single-eigenvalue) FCI problem~\cite{Cleland2010, Petruzielo2012, Greene2019, Greene2020},
and they were recently extended to the excited-state (i.e., multiple-eigenvalue) FCI problem~\cite{greene2022full}.
For example, factorizing the matrix $\mathbf{A}$ and employing additional compression operations can further reduce the cost of matrix multiplication in each iteration~\cite{Greene2019}.
Additionally, incorporating random sparsification techniques into other iterative linear algebra methods (e.g. Jacobi--Davidson \cite{Davidson1975, Sleijpen1996}) may be worth investigating, since their deterministic versions generally converge faster than subspace iteration~\cite{Blunt2019}.

Finally, despite their success in high-dimensional applications, randomized iterative algorithm are lacking explicit \emph{a priori} error bounds.
Analyzing the complicated correlations between iterates to precisely characterize their convergence properties is a pressing and ambitious goal.

\appendix

\section{Compression Algorithms}
\label{sec:compAlg}
This section provides pseudocode for implementing the pivotal compression scheme introduced in Section \ref{sec:comp} in serial and in parallel.
We also present numerical results demonstrating that pivotal compression yields less statistical error than systematic and multinomial compression schemes.

\subsection{Pivotal Compression}
\label{sec:pivComp}
As introduced in Section \ref{sec:comp}, stochastically compressing a vector $\mathbf{x}$ to $m$ nonzero elements involves the three steps that are summarized in Algorithm \ref{alg:pivComp}: identifying the elements to preserve exactly in the compressed vector, sampling from among the remaining elements, and constructing the compressed vector.
Although any unbiased sampling scheme could be used for the stochastic component of this algorithm, we specifically emphasize the pivotal sampling approach described in Algorithm \ref{alg:pivSamp}.

\begin{algorithm}[ht]
\caption{Stochastic compression with pivotal sampling}
\label{alg:pivComp}
\begin{algorithmic}[1]
\State{\textbf{Input: } A vector $\mathbf{x} \in \mathbb{R}^n$, a target number of nonzero elements $m \leq n$}
\State{Set $\mathcal{D} = \emptyset$, $d = 0$.}
\While{$\max_{i\colon i \notin \mathcal{D}} |\mathbf{x}_i| \geq \frac{1}{m-d} \sum_{j \notin \mathcal{D}} |\mathbf{x}_j|$}
\State{Set $j = \text{arg max}_{i\colon i\notin \mathcal{D}} |\mathbf{x}_i|$.}
\State{Add index $j$ to $\mathcal{D}$ and set $d = d + 1$.}
\EndWhile
\State{Set $\mathbf{p}_i = 0$ for all $i \in \mathcal{D}$ and set $\mathbf{p}_i = (m - d) |\mathbf{x}_i| \slash \sum_{j \notin \mathcal{D}} |\mathbf{x}_j|$ for all $i \notin \mathcal{D}$.}
\State{Apply pivotal sampling (Algorithm \ref{alg:pivSamp}) to sample $(m-d)$ elements $\mathcal{S}$ from $\mathbf{p}.$}
\State{Set $\Phi(\mathbf{x})_i = \mathbf{x}_i$ for all $i \in \mathcal{D}$, 
set $\Phi(\mathbf{x})_i = \mathbf{x}_i \slash \mathbf{p}_i$ for all $i \in \mathcal{S}$, and set $\Phi(\mathbf{x})_i = 0$ for all $i \notin \mathcal{D} \cup \mathcal{S}$.}
\State{\textbf{Return: } Compressed vector $\Phi(\mathbf{x})$}
\end{algorithmic}
\end{algorithm}

\begin{algorithm}[ht]
\caption{Pivotal sampling~\cite{Chauvet2017}}
\label{alg:pivSamp}
\begin{algorithmic}[1]
\State{\textbf{Input: } A number $g$ of elements to sample, a vector $\mathbf{p} \in \mathbb{R}^n$ of probabilities, with $\sum_i \mathbf{p}_i = g$}
\State{\textbf{Initialization: } Set $\mathcal{S} = \emptyset$, $b = 0$, $l = 0$, $f = 1$}
\For{$j=1,2,..., g$}
\State{Set $s = \max \{ k : b + \sum_{i=f}^k \mathbf{p}_m < 1 \}$}
\State{Randomly select index $h$ from $(l, f, f+1, f+2, ..., s)$ with probabilities proportional to $(b, \mathbf{p}_f, \mathbf{p}_{f+1}, \mathbf{p}_{f+2}, ..., \mathbf{p}_s)$}
\State{Set $a = 1 - b - \sum_{i=f}^s \mathbf{p}_i$}
\State{Set $b = \mathbf{p}_{s+1} - a$}
\State{With probability $(1 - a(1 - b)^{-1})$, add index $h$ to $\mathcal{S}$ and set $l = s+1$; otherwise, add index $(s + 1)$ to $\mathcal{S}$ and set $l = h$}
\State{Set $f = s+2$}
\EndFor
\State{\textbf{Return: } Sampled indices $\mathcal{S}$}
\end{algorithmic}
\end{algorithm}

\subsection{Parallelizing vector compression}
\label{sec:parComp}
Here we describe possible strategies for parallelizing each of the two steps involved in pivotal compression of a vector $\mathbf{x}$.
We assume that elements of $\mathbf{x}$ are distributed among $n_\text{procs}$ parallel processes, not necessarily uniformly, in arbitrary order.

Algorithm \ref{alg:pivComp} describes a serial implementation of the first step, namely identifying the largest-magnitude elements in $\mathbf{x}$ and determining the number $d$ to preserve exactly.
This procedure can be parallelized by noting that elements of $\mathbf{x}$ need not be considered in strict order of decreasing magnitude.
We provide pseudocode for this alternative approach in Algorithm \ref{alg:keepExact}, and a complete parallel implementation
is included in the open-source FRIES software on GitHub \cite{fries}.
The resulting set of deterministic indices $\mathcal{D} = \mathcal{D}^{(1)} \cup \mathcal{D}^{(2)} \cup ... \cup \mathcal{D}^{(n_\text{procs})}$, containing a total of $d$ indices, still satisfies the criteria in \eqref{eq:detCriterion} and \eqref{eq:detCriterion2}.

\begin{algorithm}
\caption{Parallel selection of elements for exact preservation}
\label{alg:keepExact}
\begin{algorithmic}[1]
\State{\textbf{Input:} Vector $\mathbf{x}^{(j)}$ on each process $j$ containing elements of the vector to be compressed, target number $m$ nonzero elements across all processes}
\State{Set $d =  0$, $\mathcal{D}^{(j)} = \emptyset$ for all $j$.} \Comment{Deterministic indices}
\For{$j=1,2,...,n_\text{procs}$} \Comment{In parallel}
\State{Set $w^{(j)} = \sum_i |\mathbf{x}^{(j)}_i|$.}
\EndFor
\State{Communicate values of $w^{(j)}$ among all processes.} \label{alg:keep-restart}
\For{$j=1,2,...,n_\text{procs}$} \Comment{In parallel}
\State{Set $a^{(j)} = 0$.}
\While{$\max_{i \colon i \notin \mathcal{D}} |\mathbf{x}^{(j)}_i| \geq \frac{1}{(m - d - a^{(j)})}\sum_{l=1}^{n_\text{procs}} w^{(l)}$}
\State{Set $l = \text{arg max}_{i \colon i \notin \mathcal{D}^{(j)}} |\mathbf{x}^{(j)}_i|$}
\State{Add index $l$ to $\mathcal{D}^{(j)}$ and set $a^{(j)} = a^{(j)} + 1$.}
\State{Set $w^{(j)} = w^{(j)} - |\mathbf{x}^{(j)}_l|$.}
\EndWhile
\EndFor
\State{Set $d = d + \sum_{j=1}^{n_\text{procs}} a^{(j)}$.}
\If{$\max_j a^{(j)} >0$}
\State{Goto line \ref{alg:keep-restart}.}
\EndIf
\State{\textbf{Return: } Deterministic indices $\mathcal{D}^{(j)}$ for each process}
\end{algorithmic}
\end{algorithm}

Next we describe an approach to parallelizing the second sampling step in vector compression.
The vector $\mathbf{p}$ of probabilities \eqref{eq:probabilities} is divided into vectors $\mathbf{q}^{(1)}$, $\mathbf{q}^{(2)}$, 
..., 
$\mathbf{q}^{(n_\text{procs})}$, 
where $\mathbf{q}^{(j)}$ contains the elements of the vector on a particular process $j$.
The number $\mathbf{g}_j$ of samples assigned to each process is a random number satisfying
\begin{equation}
\label{eq:start}
\textup{E} [\mathbf{g}_j] = (m-d) \frac{\left\lVert \mathbf{q}^{(j)} \right\rVert_1}{\left\lVert \mathbf{p} \right\rVert_1},
\end{equation}
with the additional constraint $\sum_j \mathbf{g}_j = m-d$.
We construct $\mathbf{g}$ by sampling $(m - d - \sum_j \lfloor \textup{E}[\mathbf{g}_j] \rfloor)$ indices $j$ with selection probabilities
%
\begin{equation}
\textup{P}\left\{\textup{select } j\right\} =
\textup{E}[\mathbf{g}_j] - \lfloor \textup{E}[\mathbf{g}_j] \rfloor
\end{equation}
by using pivotal sampling (Algorithm \ref{alg:pivSamp}).
Denoting the sampled indices as $\mathcal{S}'$, elements of $\mathbf{g}$ are given as
\begin{equation}
\label{eq:end}
    \mathbf{g}_j = \begin{cases}
    \lfloor \textup{E}[\mathbf{g}_j] \rfloor + 1, & j \in \mathcal{S}' \\
    \lfloor \textup{E}[\mathbf{g}_j] \rfloor, & j \notin \mathcal{S}'
    \end{cases}
\end{equation}
This operation can be performed efficiently without parallelization, assuming that the dimension of $\mathbf{t}$ (i.e.~the number of parallel processes) is small.

\begin{algorithm}
\caption{Parallel sampling of vector elements}
\label{alg:parSamp}
\begin{algorithmic}[1]
\State{\textbf{Input: } Vector $\mathbf{q}^{(j)}$ containing selection probabilities for elements on each process $j$, total number $g$ of elements to sample from all processes}
\State{Set $c = g$}
\For{$j=1,2,...,n_\text{procs}$} \Comment{$n_\text{procs}$ is the number of processes}
\State{Set $\mathbf{a}_j = g ||\mathbf{q}^{(j)}||_1 \left( \sum_i ||\mathbf{q}^{(i)}||_1 \right)^{-1}$}
\State{Set $\mathbf{t}_j = \mathbf{a}_j - \lfloor \mathbf{a}_j \rfloor$}
\State{Set $\mathbf{g}_j = \lfloor \mathbf{a}_j \rfloor$}
\State{Set $c = c - \mathbf{g}_j$}
\EndFor
\State{Sample $c$ elements $\mathcal{S}'$ from $\mathbf{t}$ by pivotal sampling (Algorithm \ref{alg:pivSamp})}
\State{Add 1 to $\mathbf{g}_j$ for each $j \in \mathcal{S}$}
\For{$j =1,2,...,n_\text{procs}$} \Comment{In parallel}
\State{Set $\mathbf{s}_j = \sum_i \mathbf{q}^{(j)}_i$}
\If{$\mathbf{g}_j > \mathbf{a}_j$}
\For{$i=1,2,...$}
\State{$\mathbf{y}_i^{(j)} = \min \left\{ 1, \mathbf{q}^{(j)}_i / \mathbf{t}_j \right\}$}
\State{Set $\mathbf{s}_j = \mathbf{s}_j + \mathbf{y}_i^{(j)} - \mathbf{q}^{(j)}_i$}
\State{Set $\mathbf{q}^{(j)}_i = \mathbf{y}_i^{(j)}$}
\If{$\mathbf{s}_j \geq \mathbf{g}_j$}
\State{Set $\mathbf{q}^{(j)}_i = \mathbf{y}_i^{(j)} + \mathbf{g}_j - \mathbf{s}_j$}
\State{Terminate for loop}
\EndIf
\EndFor
\Else
\For{$i=1,2,...$}
\State{$\mathbf{y}_i^{(j)} = \max \{0, (\mathbf{q}^{(j)}_i - \mathbf{t}_j) / (1 - \mathbf{t}_j) \}$}
\State{Set $\mathbf{s}_j = \mathbf{s}_j + \mathbf{y}_i^{(j)} - \mathbf{q}^{(j)}_i$}
\State{Set $\mathbf{q}^{(j)}_i = \mathbf{y}_i^{(j)}$}
\If{$\mathbf{s}_j \leq \mathbf{g}_j$}
\State{Set $\mathbf{q}^{(j)}_i = \mathbf{y}_i^{(j)} + \mathbf{g}_j - \mathbf{s}_j$}
\State{Terminate for loop}
\EndIf
\EndFor
\EndIf
\EndFor
\State{Sample $\mathbf{g}_j$ elements from $\mathbf{q}^{(j)}$ on each process $j$ by pivotal sampling (Algorithm \ref{alg:pivSamp}, in parallel)}
\end{algorithmic}
\end{algorithm}

Next, the probabilities on each process $j$ must be adjusted to ensure that their sum is $\mathbf{g}_j$.
This adjustment is performed differently depending on whether $\mathbf{g}_j > \textup{E}[\mathbf{g}_j]$ or $\mathbf{g}_j < \textup{E}[\mathbf{g}_j]$.
No adjustment is needed if $\mathbf{g}_j = \textup{E}[\mathbf{g}_j]$, i.e.~if $\textup{E}[\mathbf{g}_j]$ is integer-valued.
We define the vectors $\mathbf{y}^{(j)}$, with elements
\begin{equation}
\mathbf{y}^{(j)}_i = \begin{cases}
\min \left\lbrace 1, \mathbf{q}^{(j)}_i / t_j \right\rbrace & \mathbf{g}_j > \textup{E}[\mathbf{g}_j] \\
\max \left\lbrace 0, \left(\mathbf{q}^{(j)}_i - \mathbf{t}_j \right) / \left(1 - \mathbf{t}_j \right) \right\rbrace & \mathbf{g}_j < \textup{E}[\mathbf{g}_j]
\end{cases}
\end{equation}
and $\mathbf{z}^{(j)}$, with elements
\begin{equation}
\mathbf{z}^{(j)}_i = \sum_{l=1}^i \mathbf{y}^{(j)}_l + \sum_{l = i + 1}^{{e}_j} \mathbf{q}^{(j)}_l
\end{equation}
where $e_j$ is the total number of elements in the vector $\mathbf{q}^{(j)}$.
The index $h_j$ is calculated as the minimum value of $i$ satisfying $\mathbf{z}^{(j)}_i \geq \mathbf{g}_j$ if $\mathbf{g}_j > \textup{E}[\mathbf{g}_j]$ or $\mathbf{z}^{(j)}_i \leq \mathbf{g}_j$ if $\mathbf{g}_j < \textup{E}[\mathbf{g}_j]$.
The adjusted probabilities are then calculated as
\begin{equation}
\mathbf{q}^{(j)\prime}_i = 
\begin{cases}
\mathbf{y}^{(j)}_i & i < h_j \\
\mathbf{g}_j - \sum_{l=1}^{h_j} \mathbf{y}^{(j)}_l - \sum_{l=h+1}^{e_k} \mathbf{q}^{(j)}_l & i = h_j \\
\mathbf{q}^{(j)}_i & i > h_j
\end{cases}
\end{equation}
This particular approach was chosen to minimize the number of probabilities to be recalculated.

Finally, after calculating the adjusted probabilities, we sample $\mathbf{g}_j$ elements according to the pivotal scheme in Algorithm \ref{alg:pivSamp}.

\subsection{Comparing compression schemes}
\label{sec:pivVsys}
Here we compare pivotal compression to two other compression schemes we have investigated previously, namely multinomial and systematic compression~\cite{Greene2019}.
In the multinomial scheme, no elements are preserved exactly: instead, the compressed vector is constructed by independently sampling indices from the input vector with weights proportional to the magnitudes of the corresponding elements.
In systematic compression, some elements are preserved exactly according to \eqref{eq:detCriterion} and \eqref{eq:detCriterion2}, and the remaining elements are sampled using a systematic sampling approach.

\begin{figure}
    \centering
    \includegraphics[scale=0.75]{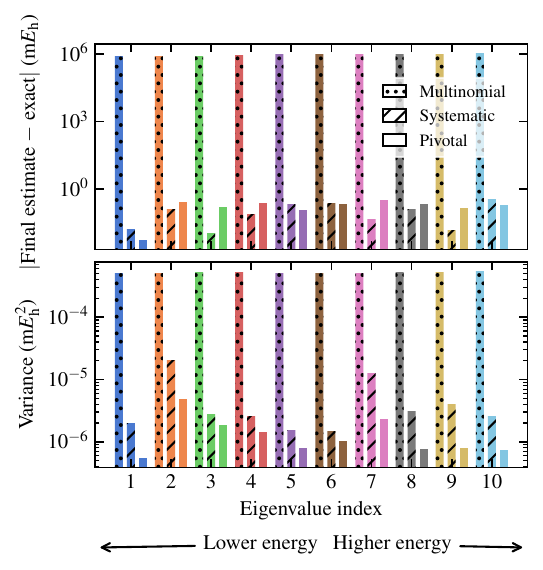}
    \caption{Results obtained by applying our randomized subspace iteration with three different compression schemes (multinomial, systematic, and pivotal) to estimate eigenvalues for the Ne atom. In all schemes, iterate matrices were compressed to $m=10,000$ nonzero elements per column. (top) The magnitude of the bias for each estimate, obtained from a trajectory of 50,000 iterations (with $i_\text{min} = 20,000$). (bottom) The variance for each estimate, obtained using \eqref{eq:varEstimate}. Note the logarithmic scale on the vertical axes.}
    \label{fig:multisyspivVariance}
\end{figure}

Figure \ref{fig:multisyspivVariance} presents results obtained by applying our randomized subspace iteration with each of these compression schemes to the Ne system described above.
All estimates were obtained from trajectories of 50,000 iterations with the burn-in time $i_\text{min}$ chosen as 20,000.
Biases from randomized subspace iteration with multinomial compression are approximately 7 orders of magnitude greater than with systematic or pivotal schemes.
The variance resulting from multinomial compression is approximately 50 times greater the variance from the systematic scheme.
Biases from the systematic and pivotal schemes were approximately the same, 
while the pivotal scheme exhibits approximately 3.5 times less variance.
These results provide further justification for our use of pivotal compression.

\section{Stability analysis of numerical experiments}
\label{sec:stability}
\begin{figure}
    \centering
    \includegraphics[width=0.5\linewidth]{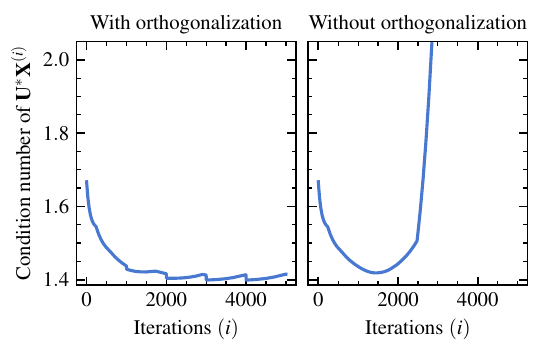}
    \caption{The condition number of the matrix $\mathbf{U}^{\ast} \mathbf{X}^{(i)}$ at each iteration $i$ of deterministic subspace iteration calculations on the Ne atom, either performing orthogonalization every 1000 iterations (left) or not at all (right).}
    \label{fig:deterministicCN}
\end{figure}

Plots of the condition number of the matrix $\mathbf{U}^{\ast} \mathbf{X}^{(i)}$ can be used to monitor the stability of calculations as the iteration proceeds.
An increasing trend in the condition number can indicate that orthogonalization is being performed too infrequently.
In order to illustrate this, we present plots of this condition number obtained by applying two versions of deterministic subspace iteration to the Ne system introduced in Section \ref{sec:detNe}.
In the calculation shown in the left panel of Figure \ref{fig:deterministicCN}, orthogonalization is performed at intervals of 1000 iterations, while in the right panel, orthogonalization is not performed at all.
Although both of these calculations yield the same eigenvalue estimates, condition numbers in the calculation without orthogonalization increase rapidly, which suggests that the algorithm would encounter numerical instabilities if run for more iterations.
After 5000 iterations, the condition number for this calculation is 59.
In contrast, condition numbers in the calculation with orthogonalization remain less than 1.5.

\begin{figure}[h]
    \centering
    \includegraphics[width=0.5\linewidth]{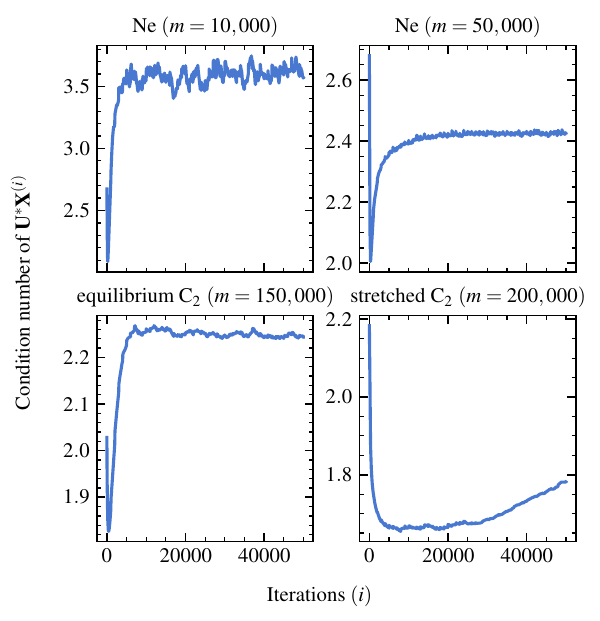}
    \caption{The condition number of $\mathbf{U}^{\ast} \mathbf{X}^{(i)}$ at each iteration $i$ for each of the randomized calculations presented in the main text.}
    \label{fig:condNum}
\end{figure}

Similarly, we monitored condition numbers for each calculation performed with our randomized subspace iteration (Sections \ref{sec:detNe} and \ref{sec:C2results}).
These are presented in Figure \ref{fig:condNum}.
In all cases except for stretched \ce{C2}, the condition number stabilized.
A gradual increasing trend was observed at later iterations for stretched \ce{C2}.
We anticipate that this condition number would stabilize after more iterations and that the increase observed here did not affect the accuracy of our results, especially considering their close agreement with the numerically exact eigenvalues.

\begin{figure}
    \centering
    \includegraphics[width=0.75\linewidth]{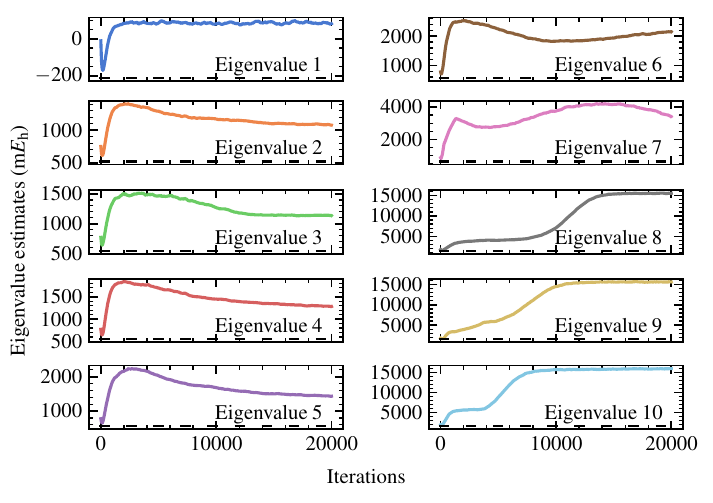}
    \caption{\textbf{Eigenvalue estimates from randomizing a standard subspace iteration.} Estimates of the ten least-energy eigenvalues for Ne were obtained by compressing iterate matrices to $m=10,000$ nonzero elements per column. Note the significantly greater errors relative to our randomized non-standard subspace iteration (Fig. \ref{fig:NeEst}). Dashed lines indicate exact eigenvalues.}
    \label{fig:symmEvals}
\end{figure}

\section{Randomization of standard subspace iteration fails}
\label{sec:symmSubsp}
This section describes the application of repeated stochastic sparsification techniques to a more standard subspace iteration that relies more on nonlinear operations on the iterates $\mathbf{X}^{(i)}$.
In the nonlinear approach, multiplication by $[\mathbf{G}^{(i)}]^{-1}$ enforces orthonormality in the full $n$-dimensional vector space.
At regular intervals, $\mathbf{G}^{(i)}$ is constructed from a Gram-Schmidt orthogonalization of $\mathbf{X}^{(i)}$.
At other iterations, $\mathbf{G}^{(i)}$ is a diagonal matrix containing the $\ell_1$-norms of the columns of $\mathbf{X}^{(i)}$.
For these tests, we performe orthogonalization at intervals of 1000 iterations.

In each iteration, eigenvalues are estimated by applying the Rayleigh-Ritz method to the compressed iterate $\Phi(\mathbf{X}^{(i)})$, i.e.~by solving the generalized eigenvalue equation
\begin{equation}
\label{eq:quadEst}
\Phi(\mathbf{X}^{(i)})^* \mathbf{A} \Phi(\mathbf{X}^{(i)}) \mathbf{W}^{(i)} = \Phi(\mathbf{X}^{(i)})^* \Phi(\mathbf{X}^{(i)}) \mathbf{W}^{(i)} \mathbf{\Lambda}^{(i)}
\end{equation}
for the matrix $\mathbf{\Lambda}^{(i)}$ of Ritz values.
Because this equation involves quadratic inner products of the compressed iterates, the resulting eigenvalue estimates are variational.
Applying this approach with iterates compressed to $m=10,000$ to the Ne system defined in the main text yields the eigenvalue estimates presented in Fig. \ref{fig:symmEvals}.
The best (i.e.~minimum) estimates differ from the exact eigenvalues by as much as 201 m$E_\text{h}$.

One might expect that averaging can be used to improve the accuracy of these eigenvalue estimates.
However, because this approach is variational, averaging the Ritz values themselves yields estimates with at least as much error as the minimum eigenvalue estimates considered above.
Instead averaging the $k \times k$ matrices $\Phi(\mathbf{X}^{(i)})^* \mathbf{A} \Phi(\mathbf{X}^{(i)})$ and $\Phi(\mathbf{X}^{(i)})^* \Phi(\mathbf{X}^{(i)})$ yields poorer estimates.
These differ from the exact eigenvalues by 299 to 10,590 m$E_\text{h}$.
Although the computational costs and memory requirements for this nonlinear approach are approximately the same as for the method in the main text, it is impossible to extract eigenvalue estimates of similar accuracy.
This underscores the importance of the non-standard choices made in constructing our randomized subspace iteration.

Previous methods have sought to reduce the variance of estimated eigenvalues in alternative ways.
For example, the bias could potentially be reduced by replacing one instance of the iterate matrix in \eqref{eq:quadEst} by an independently generated replica trajectory and omitting the orthogonalization step.
This strategy has been tested previously in the context of the FCI problem \cite{Blunt2018nonlinear}.
However, the variance in the inner products between vectors from independent trajectories can scale unfavorably with the matrix dimension $n$. 

Another variance reduction scheme, introduced by \cite{Ceperley1988}, 
is an early precursor to the algorithm called the `variational approach for conformation dynamics' (VAC) \cite{Noe2013,webber2021error,lorpaiboon2020integrated}. 
In 
VAC, one needs to select a ``time-lag'' parameter sufficiently long to avoid bias due to the imperfect choice of trial vectors, but sufficiently short to avoid collapse to the dominant eigenvector.
The difficulty in tuning the lag-time parameter is one challenge that may hinder the application of this technique. 

\section{Data Availability}
All data from the numerical experiments presented here is available at \url{https://doi.org/10.5281/zenodo.4624477}.
The code used to perform our numerical experiments can be accessed at \url{https://github.com/sgreene8/FRIES}.

\section*{Acknowledgments}
We gratefully acknowledge productive discussions with Aaron Dinner, Michael Lindsey, Verena Neufeld, Joel Tropp, Ethan Epperly, and James Smith throughout the development and execution of this project. Lek-Heng Lim originally raised the possibility of randomizing subspace iteration to us. Benjamin Pritchard provided invaluable suggestions for improving the readability and efficiency of our source code. Computational resources were provided by the Research Computing Center at the University of Chicago and the High Performance Computing Center at New York University.

\bibliographystyle{siamplain}
\bibliography{Sams_refs, software}
\end{document}